\documentclass[aap,citesort,MSNbibl,dvips]{arximspdf}
\usepackage{mathbh}
\usepackage{graphicx}

\doi{10.1214/10-AAP734}
\volume{21}
\issue{4}
\pubyear{2011}
\firstpage{1215}
\lastpage{1252}

\makeatletter

\newtheorem{theor}{Theorem}
\newtheorem{lemma}[theor]{Lemma}

\newcommand{\N}{\mathbb{N}}
\newcommand{\Z}{\mathbb{Z}}
\newcommand{\GH}{\mathcal{H}}
\newcommand{\K}{\mathbb{K}}
\newcommand{\E}{\mathbb{E}}
\newcommand{\ep}{\varepsilon}
\newcommand{\ind}{\mathbh1}
\newcommand{\norm}[1]{|\!|#1|\!|}
\newcommand{\tow}{\rightleftharpoons}
\newcommand{\ver}{\updownarrow}
\newcommand{\hor}{\leftrightarrow}
\newcommand{\W}{\mathcal{W}}
\newcommand{\X}{\mathcal{X}}
\newcommand{\Y}{\mathcal{Y}}

\newcommand{\card}{\operatorname{card}}

\def\bptnote#1{}

\makeatother

\begin{document}
\begin{frontmatter}

\title{Contact and voter processes on the infinite percolation cluster
as models of host-symbiont interactions}
\runtitle{Host-symbiont interactions}

\begin{aug}
\author[A]{\fnms{D.} \snm{Bertacchi}},
\author[B]{\fnms{N.} \snm{Lanchier}\corref{}\thanksref{t1}\ead[label=e1]{lanchier@math.asu.edu}}
and
\author[C]{\fnms{F.} \snm{Zucca}}
\runauthor{D. Bertacchi, N. Lanchier and F. Zucca}
\affiliation{Universit\`a di Milano-Bicocca, Arizona State University
and Politecnico~di~Milano}
\address[A]{D. Bertacchi\\
Dipartimento di Matematica e Applicazioni\\
Universit\`a di Milano-Bicocca \\
Via Cozzi 53\\
20125 Milano\\
Italy}
\address[B]{N. Lanchier\\
School of Mathematical\\
\quad and Statistical Sciences\\
Arizona State University\\
Tempe, Arizona 85287\\
USA\\
\printead{e1}}
\address[C]{F. Zucca\\
Dipartimento di Matematica\\
Politecnico di Milano\\
Piazza Leonardo da Vinci 32\\
20133 Milano\\
Italy}
\end{aug}

\thankstext{t1}{Supported in part by NSF Grant DMS-10-05282.}

\received{\smonth{11} \syear{2009}}
\revised{\smonth{9} \syear{2010}}

%
\begin{abstract}
We introduce spatially explicit stochastic processes to model
multispecies host-symbiont interactions.
The host environment is static, modeled by the infinite percolation
cluster of site percolation.
Symbionts evolve on the infinite cluster through contact or voter type
interactions, where each host may be infected by a
colony of symbionts.
In the presence of a single symbiont species, the condition for
invasion as a function of the density of the habitat of
hosts and the maximal size of the colonies is investigated in details.
In the presence of multiple symbiont species, it is proved that the
community of symbionts clusters in two dimensions whereas
symbiont species may coexist in higher dimensions.
\end{abstract}

\setattribute{keyword}{AMS}{AMS 2000 subject classification.}

\end{frontmatter}

\section{Introduction}
\label{sec:introduction}

The term symbiosis was coined by the mycologist Heinrich Anto de Bary
to denote close and long-term physical and biochemical
interactions between different species, in contrast with competition
and predation that imply only brief interactions.
Symbiotic relationships involve a~symbiont species, smaller in size,
that always benefits from the relationship, and a host species,
larger in size, that may either suffer, be relatively unaffected, or
also benefit from the relationship, which are referred to as
parasistism, commensalism, and mutualism, respectively.
The degree of specificity of the symbiont is another important factor:
while some symbionts may live in association with a wide range of host
species, in which case the symbiont is called a generalist,
others are highly host-specific indicating that they can only benefit
from few host species.
Symbiotic relationships, either pathogenic or mutualistic, are
ubiquitous in nature.
For instance, more than 90\% of terrestrial plants \cite{trappe1987}
live in association with mycorrhizal fungi, with the plant
providing carbon to the fungus and the fungus providing nutrients to
the plant, most herbivores have mutualistic gut fauna that
help them digest plant matter, and almost all free-living animals are
host to one or more parasite taxa \cite{price1980}.

To understand the role of spatial structure on the persistence of
host-parasite and host-mutualist associations, Lanchier
and Neuhauser \cite
{lanchierneuhauser2006a,lanchierneuhauser2006b,lanchierneuhauser2009}
have initiated the study of
multispecies host-symbiont systems including local interactions based
on interacting particle systems.
The stochastic process introduced in \cite{lanchierneuhauser2006b}
describes the competition among specialist and generalist
symbionts evolving in a deterministic static environment of hosts.
The mathematical analysis of this model showed that fine-grained
habitats promote generalist strategies, while coarse-grained
habitats increase the competitiveness of specialists.
The stochastic process introduced in \cite
{lanchierneuhauser2006a,lanchierneuhauser2009} includes in addition a feedback
of the hosts, which is modeled by a dynamic-host system.
This process has been further extended by Durrett and Lanchier~\cite
{durrettlanchier2008}.
The host population evolves, in the absence of symbionts, according to
a~biased voter model, while the symbiont population
evolves in this dynamic environment of hosts according to a contact
type process.
The parameters of the process allow to model the effect of the
symbionts on their host as well as the degree of specificity
of the symbionts, thus resulting in a system of coupled interacting
particle systems, each describing the evolution of
a~trophic level.
The model is designed for the understanding of the role of the
symbionts in the spatial structure of plant communities.
It is proved theoretically that generalist symbionts have only a
limited effect on the spatial structure of
their habitat \cite{lanchierneuhauser2009}.
In contrast, the inclusion of specialist parasites promotes coexistence
of the hosts in terms of the existence of
a~stationary distribution under which the density of each host type is
positive, while the analysis of the corresponding
mean-field model supported by numerical simulations suggests that in
any dimension the inclusion of specialist mutualists
translates into a clustering of the host environment \cite
{durrettlanchier2008}.

Similarly to most spatial epidemic models such as the contact process,
the state space of the stochastic processes
introduced in \cite
{durrettlanchier2008,lanchierneuhauser2006a,lanchierneuhauser2006b,lanchierneuhauser2009}
indicates
whether hosts are either healthy or infected, but does not distinguish
between different levels of infection of the hosts.
However, it is known from past research that the number of symbiont
individuals, including ectosymbionts, that is, symbionts living
on their hosts or in their skin, associated to a single host individual
may vary significantly.
Mooring and Samuel \cite{mooringsamuel1998} found for instance an
average of 1791 individuals of the species \textit{Dermacentor
albipictus}, commonly known as Winter Tick, on individual elk in
Alberta, while some individual moose have been found with more
than 50,000 ticks.
In addition, symbionts are generally much smaller organisms than their
hosts and reproduce much faster and in greater number.
This motivates the development of spatially explicit multiscale models
of host-symbiont interactions that describe the
presence of symbionts through a level of infection of the hosts rather
than binary random variables (infected versus
healthy hosts) and include both inter-host symbiont dynamics and
intra-host symbiont dynamics.

In diversity ecology, the infrapopulation refers to all the parasites
of one species in a single individual host, while
the metapopulation refers to all the parasites of one species in the
host population.
In systems involving multiple species of parasites, all the parasites
of all species in a single individual host and in an
entire host population are called infracommunity and component
community, respectively.
This terminology shall be employed in this article for symbionts in
general, that is parasites, commensalists and mutualists,
even though, strictly speaking, it only applies to parasites.
Our main objective is to deduce from the microscopic evolution rules of
the symbionts, described by transmission rates and
reproduction rates, the long-term behavior of the metapopulation in a
single-species invasion model, and the long-term behavior
of the component community in a multispecies competition model.
Since a host species and a symbiont species involved in a symbiotic
relationship usually evolve at very different time
scales (symbionts reproduce much faster than their hosts), we shall
assume in both invasion and competition models that the
discrete habitat of hosts is static.
This habitat will be modeled by a realization of the infinite
percolation cluster of supercritical site
percolation \cite{grimmett1989}.
We shall also assume that symbionts can only survive when associated
with a host (obligate relationship), which restricts
their habitat to the infinite percolation cluster, and, to understand
the role of space on the persistence of the symbiotic
relationship, that symbionts can only transmit to nearby hosts, adding
to the complexity of the interactions.
In the single-species model, infrapopulations will evolve according to
the logistic growth process, and the entire
metapopulation according to a~mixture of this model and its spatial
analog, the contact process \cite{harris1974}.
In the multispecies model, we will assume that infracommunities evolve
according to the Moran
model~\cite{moran1958}, and the entire component community according
to a mixture of this model and its spatial analog, the
voter model \cite{cliffordsudbury1973,holleyliggett1975}.
Our analysis shows that the condition for survival of a metapopulation
strongly depends on the carrying capacity of each
infrapopulation.
Exact calculations of the critical curve as a function of the
reproduction and transmission rates are given when
infrapopulations can be arbitrarily large which, as mentioned above, is
a realistic biological assumption in many symbiotic
relationships.
In systems involving multiple symbiont species, the long-term behavior
of the component community depends on the spatial
dimension: the community clusters in two dimensions whereas coexistence
is possible in higher dimensions.


\section{Models and results}
\label{sec:model}

The models are constructed in two steps.
First, the static random environment of hosts is fixed from a
realization of the infinite percolation cluster of site
percolation \cite{grimmett1989}.
This random environment naturally induces a random graph.
The symbionts are then introduced into this universe where they evolve
according to an interacting particle system on the random
graph.
The interactions are modeled based on two of the simplest particle
systems: the contact process \cite{harris1974}
and the voter model \cite{cliffordsudbury1973,holleyliggett1975}.
The structure of the random graph implies that the infrapopulation
dynamics are described by logistic growth processes, that is,
contact processes on a~complete graph, and the infracommunity dynamics
by Moran models, that is, voter models
on a complete graph.


\subsection*{Host environment}
To define the habitat of hosts, we set $p \in(0, 1]$ and let $\omega$
be a realization of the site percolation
process with parameter $p$ on the $d$-dimensional regular lattice $\Z
^d$, that is, each site of the lattice is either permanently
occupied by an individual host with probability $p$ or permanently
empty with probability $1 - p$.
Let $\mathbb{H}(\omega)$ denote the set of open/occupied sites.
By convention, elements of $\Z^d$ and processes with state space $S
\subset\Z^d$ will be denoted in the following by capital
Latin letters.
We say that there is an open path between site $X$ and site $Y$ if
there exists a sequence of sites
$X = X_0, X_1, \ldots, X_n = Y$ such that the following two conditions hold:
\begin{enumerate}
\item For $i = 0, 1, \ldots, n$, we have $X_i \in\mathbb{H}(\omega)$, that
is, site $X_i$ is open.
\item For $i = 0, 1, \ldots, n - 1$, we have $X_i \sim X_{i + 1}$,
\end{enumerate}
where $X_i \sim X_{i + 1}$ means that the Euclidean norm $\norm{X_i -
X_{i + 1}} = 1$.
Writing $X \tow Y$ the event that sites $X$ and $Y$ are connected by an
open path, we observe that the binary relation $\tow$
is an equivalence relation on the random set $\mathbb{H}(\omega)$ thus
inducing a partition of $\mathbb{H}(\omega)$.
In dimensions $d \geq2$, there exists a critical value $p_c \in(0,
1)$ that depends on $d$ such that if $p > p_c$ then $\mathbb{H}(\omega)$
contains a unique infinite open cluster.
The infinite open cluster is also called infinite percolation cluster
and is denoted by $C_{\infty} (\omega)$ later.
We assume that $p > p_c$ from now on.
Sometimes, the infinite percolation cluster will be identified with the
graph with vertex set $C_{\infty} (\omega)$ obtained by
drawing an edge between sites of the cluster at Euclidean distance~1
from each other.
For more details about site percolation, we refer the reader to
Grimmett~\cite{grimmett1989}.


\subsection*{Random graph structure}
In order to define the state space and dynamics of the stochastic
processes, we first define a random graph $\GH(\omega)$ as follows.
Vertices of $\GH(\omega)$ are to be interpreted as possible locations
for the symbionts, while edges indicate how symbionts interact.
Let $N$ be an integer and $\K_N = \{1, 2, \ldots, N \}$.
The vertex set of $\GH(\omega)$ is
\[
C_N (\omega) = \{(X, i) \dvtx X \in C_{\infty} (\omega) \mbox
{ and } i \in\K_N \}.
\]
By convention, elements of and processes with state space $C_N (\omega
)$ will be denoted by small Latin letters.
To define the edge set, we also introduce
\[
\pi\dvtx C_N (\omega) \longrightarrow C_{\infty} (\omega)
\!\qquad\mbox{defined by } \pi(x) = X \!\qquad\mbox{for all } x = (X, i)
\in C_N (\omega).
\]
That is, $\pi(x)$ is the $C_{\infty} (\omega)$-coordinate of vertex $x$.
Let $x, y \in C_N (\omega)$.
Then vertices $x$ and $y$ are connected by an edge if and only if one
of the following two cases occurs:
\begin{enumerate}
\item If $\pi(x) = \pi(y)$, then $x$ and $y$ are connected by a
vertical edge: we write $x \ver y$.
It is convenient to assume that each vertex is connected to itself by a
vertical edge.
\item If $\pi(x) \sim\pi(y)$, then $x$ and $y$ are connected by a
horizontal edge: we write $x \hor y$.
\end{enumerate}
In words, a complete graph with $N$ vertices (which are connected to
themselves) is attached to each site of the infinite
percolation cluster.
Edges of these complete graphs are said to be vertical while, for any
two sites of the infinite percolation cluster, vertices of
the corresponding complete graphs are connected by edges which are said
to be horizontal.
Vertical and horizontal edges correspond, respectively, to potential
reproduction events and transmission events of the symbionts.



\subsection*{Invasion of a single symbiont---contact process}
\begin{figure}
\begin{tabular}{@{}cc@{}}

\includegraphics{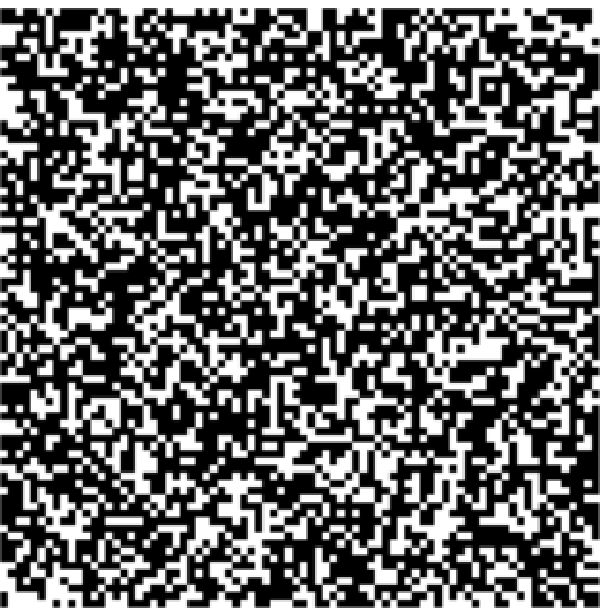}
 & \includegraphics{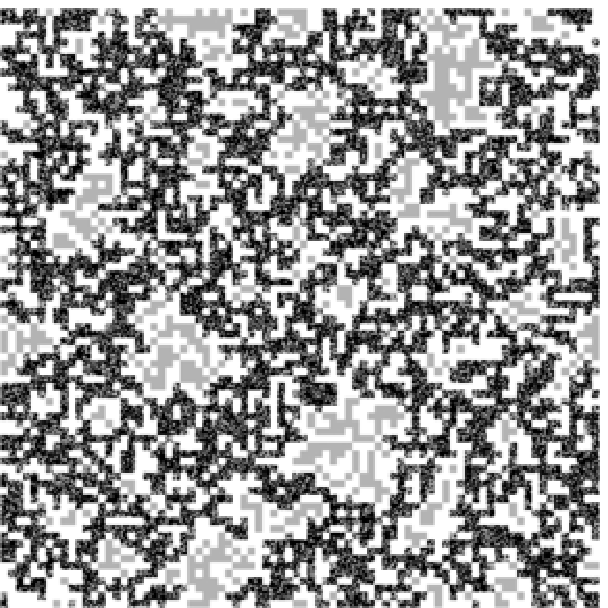}\\
(a) & (b)
\end{tabular}%
\caption{\textup{(a)} Realization of site percolation with parameter
$p = 0.6$ on the $80 \times80$ torus, with black squares referring
to open sites, and white squares to closed sites.
\textup{(b)} Snapshot of the invasion model on the percolation
structure starting with a single infected host at the center of the universe
($\alpha= 1$ and $\beta= 2$).
Each site is represented by a $5 \times5$ square, that is, complete graph
with $N = 25$ vertices.
White squares refer to empty sites, that is, sites which are not occupied
by a host, black dots refer to symbionts, and gray dots
to empty vertices.}
\label{fig:SHMC}
\end{figure}

To understand the conditions for survival of a single symbiont species,
we introduce a generalization of the
contact process \cite{harris1974} on the infinite random graph $\GH
(\omega)$.
This defines a continuous-time Markov process whose state space
consists of the set of the spatial configurations
$\eta\dvtx C_N (\omega) \longrightarrow\{0, 1 \}$, and whose dynamics
are described by the Markov generator $L_1$ defined on
the set of the cylinder functions by
\begin{eqnarray*}
L_1 f (\eta) &=& \sum_{x \in C_N (\omega)} [f
(\eta_{x, 0}) - f (\eta)] \\
&&{} +
\sum_{x \in C_N (\omega)} \biggl(\frac{\alpha}{N}
\sum_{x \ver y} \eta(y) +
\frac{\beta}{N \deg\pi(x)} \sum_{x \hor y} \eta(y) \biggr)\\
&&\hspace*{45.2pt}{}\times
[f (\eta_{x, 1}) - f (\eta)],
\end{eqnarray*}
where $\deg\pi(x)$ is the degree of $\pi(x)$ as a site of the
cluster $C_{\infty} (\omega)$, and where $\eta_{x, i}$ is the
configuration obtained from $\eta$ by assigning the value $i$ to
vertex $x$.
Note that the degree of each site of the infinite percolation cluster
is at least~1, therefore the dynamics are well defined.
Thinking of vertices in state 0 as uninfected and vertices in state 1
as infected by a symbiont, the expression of the
Markov generator above indicates that symbionts die independently of
each other at rate 1, reproduce within their host at the
reproduction rate $\alpha$, and transmit their offspring to the nearby
hosts at the transmission rate $\beta$.
That is, each symbiont gives birth at rate $\alpha$ to an offspring
which is then sent to a vertex chosen uniformly at random
from the parent's host.
If the vertex is uninfected, then it becomes infected while if it is
already infected then the birth is suppressed.
Similarly, each symbiont gives birth at rate $\beta$ to an offspring
which is then sent to a vertex chosen uniformly at random
from the hosts adjacent to the parent's host, which results as
previously in an additional infection if and only if the vertex
is not already infected.
See Figure \ref{fig:SHMC}
for simulation pictures of this contact process.
To study the single-species model, we will sometimes consider the
stochastic process
%
%
\[
\bar\eta_t (X) = \sum_{\pi(x) = X} \eta_t (x)\qquad \mbox
{for all } X \in C_{\infty} (\omega),
\]
where the sum is over the vertices $x \in C_N (\omega)$ such that $\pi
(x) = X$.
That is, $\bar\eta_t (X)$ keeps track of the level of infection of
the host at $X$.
This defines a~Markov process whose state space consists of the
functions that map~$C_{\infty} (\omega)$
into $\{0, 1, \ldots, N \}$ and whose dynamics are described by
\begin{eqnarray*}
\bar L_1 f (\bar\eta) &=& \sum_{X \in C_{\infty}
(\omega)} \bar\eta(X) [f (\bar\eta_{X-}) - f (\bar\eta)] \\
&&{} +
\sum_{X \in C_{\infty} (\omega)} \biggl(1 - \frac
{\bar\eta(X)}{N} \biggr) \biggl(\alpha\bar\eta(X) +
\frac{\beta}{\deg(X)} \sum_{X \sim Y} \bar\eta(Y) \biggr)\\
&&\hspace*{48.3pt}{}\times [f
(\bar\eta_{X+}) - f (\bar\eta)],
\end{eqnarray*}
where the configurations $\bar\eta_{X-}$ and $\bar\eta_{X+}$ are
obtained from the configuration~$\bar\eta$ by, respectively,
removing and adding a symbiont at site $X$.
In view of the geometry of the graph $\GH(\omega)$, the stochastic
process $\{\bar\eta_t \}_t$ can be seen as a~mixture of the contact
process with infection parameter $\beta$ on the infinite percolation
cluster and logistic growth processes with parameter $\alpha$.

To describe the predictions based on the invasion model, we let $\delta
_i$ be the measure that concentrates on
the ``all $i$'' configuration restricted to $C_N (\omega)$, that is,
\[
\delta_i \{\eta(x) = i \} = 1 \qquad\mbox{for all } x \in C_N
(\omega) = C_{\infty} (\omega) \times\K_N.
\]
We denote by $\bar\mu$ the upper invariant measure of the process $\{
\eta_t \}_t$, which is also the limit starting from the measure
$\delta_1$ since the process is attractive.
The process or metapopulation is said to survive whenever $\bar\mu
\neq\delta_0$ and is said to die out otherwise.

First, we observe that, starting with a single infection at time 0, the
number of symbionts in the system is dominated
stochastically by the number of individuals in a birth and death
process with birth parameter $b = \alpha+ \beta$ and death
parameter~1.
Recurrence of one-dimensional symmetric random walks implies that such
a process eventually dies out when $b \leq1$.
It follows that $\{\eta_t \}_t$ dies out for all values of $N$
whenever $\alpha+ \beta\leq1$.

To find a general condition for survival of the infection, we now
assume that $N = 1$ so that the value of the reproduction rate
$\alpha$ becomes irrelevant, and compare the process with the
one-dimensional contact process.
Let $\Gamma$ be an arbitrary infinite self-avoiding path in the
infinite percolation cluster~$C_{\infty} (\omega)$.
Since for all sites $X \in\Gamma$ the degree of $X$ ranges from 2 to
$2d$, the process restricted to the infinite path $\Gamma$, that is, symbionts
sent outside $\Gamma$ are instantaneously killed, dominates
stochastically the contact process on $\Gamma$ with infection parameter
$\beta/ d$.
It follows that the process survives whenever $\beta> d \beta_c
(1)$ where $\beta_c (1)$ is the critical value of the one-dimensional
contact process, since the self-avoiding path is isomorphic to~$\Z$.
Standard coupling arguments also imply that the survival probability of
the infection is nondecreasing with respect to both the
reproduction rate $\alpha$ and the maximum number of symbionts per
host $N$.
It follows directly from these monotonicity properties that, for all
values of $N$ and $\alpha$, survival occurs
whenever $\beta> d \beta_c (1)$.

We now look at the long-term behavior of the metapopulation when $N$ is large.
As previously explained, this assumption is realistic in a number of
symbiotic relationships, including the interactions between
moose and Winter Ticks \cite{mooringsamuel1998}.
Under this assumption, at least when the number of symbionts is not too
large, the stochastic process looks locally like a branching
random walk on the random graph~$\GH(\omega)$, namely the process
modified so that births onto infected vertices are allowed.
In the context of large infrapopulations, global survival of the
metapopulation occurs when the reproduction rate $\alpha> 1$ and
the transmission rate $\beta> 0$.
This and the comparison with a~birth and death process imply that, when
$N$ is large and the transmission rate $\beta$
is small, a situation which is common in parasitic relationships, the
metapopulation undergoes a phase transition when the
reproduction rate~$\alpha$ approaches 1.
Provided the density of the habitat is large enough, the phase
transition occurs more generally when the sum of the reproduction and
transmission rates approaches 1.
These results are summarized in the following theorem where
``survival'' means strong survival of the stochastic process, that is
the existence of a stationary distribution under which the density of
symbionts is positive.
\begin{theor}[(Contact interactions)]
\label{invasion}
Assume that $p > p_c$ and $\beta> 0$.
\begin{enumerate}
\item For all $N > 0$, the metapopulation dies out if $\alpha+ \beta
\leq1$ while it survives if $\beta/ d > \beta_c (1)$.
%
\item If $\alpha+ \beta> 1$ and $p$ is close to 1, then the
metapopulation survives for $N$ large.
\item If $\alpha+ \beta/ d > 1$ and $p > p_c$, then the
metapopulation survives for $N$ large.
\end{enumerate}
\end{theor}

As previously explained, the first statement of part 1 follows from a
comparison with a two-parameter branching random walk, and the second
statement from a comparison with the contact process restricted to a
self-avoiding path embedded in the infinite percolation cluster.
The proof of the second part relies on the combination of random walk
estimates and block constructions to compare the process view under suitable
space and time scales with oriented percolation, and we refer to
Section \ref{sec:invasion} for more details.
Survival when $\alpha> 1$ and $p > p_c$ in the presence of large
infrapopulations can be proved based on estimates for the extinction time
of the logistic growth process and a new block construction.
However, the third part, which indicates survival under the weaker
assumption $\alpha+ \beta/ d > 1$, can be directly deduced from the
proof of
the second part by again looking at the process restricted to an
infinite self-avoiding path of hosts.
Let $\Gamma$ be an infinite self-avoiding path, which exists almost
surely under the assumption $p > p_c$, and observe that, since the
degree of each site along this path ranges from 2 to $2d$, the process
restricted to $\Gamma$ dominates stochastically the one-dimensional
process with parameters $\alpha$ and $\beta/ d$.
The latter survives if $\alpha+ \beta/ d > 1$ since, under this
assumption, the proof of the second part indicates that survival occurs when
$p = 1$ in any dimension, including $d = 1$.
The third part of the theorem clearly follows.





\subsection*{Competition among multiple symbionts---voter model}
\begin{figure}
\begin{tabular}{@{}cc@{}}

\includegraphics{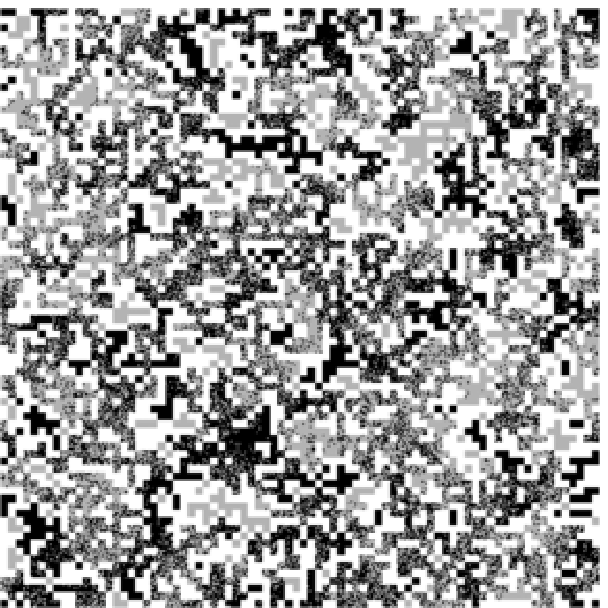}
 & \includegraphics{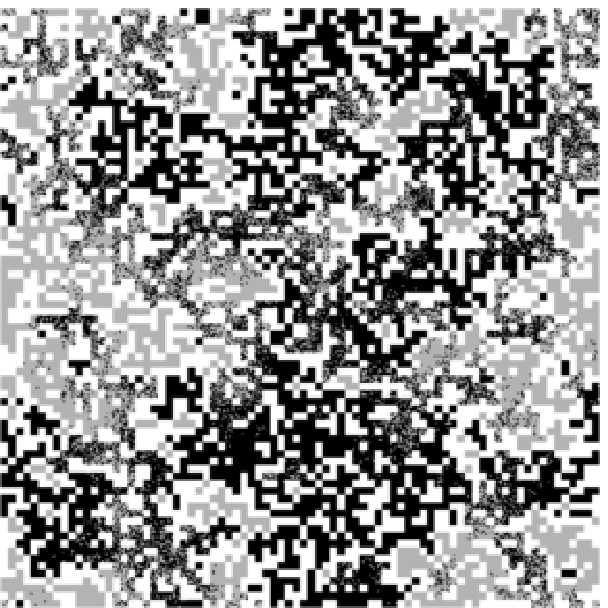}\\
(a) & (b)
\end{tabular}
\caption{Snapshots at time 100 \textup{(a)} and at time 1000 \textup
{(b)}, respectively, of the
neutral competition model on site percolation with parameter
$p = 0.6$ starting from a Bernoulli product measure with density $1/2$
($\alpha_1 = \alpha_2 = 1/2$ and $\beta_1 = \beta_2
= 1/2$).
Each site of the lattice is represented by a $5 \times5$ square, that is,
complete graph with $N = 25$ vertices.
White squares refer to empty sites, that is, sites which are not occupied
by a host, and black and gray dots to symbionts of type 1
and 2, respectively.}
\label{fig:SHMV}
\end{figure}

To study the interactions among multiple symbiont species, we introduce
the analog of the previous model replacing contact
interactions with voter interactions \cite
{cliffordsudbury1973,holleyliggett1975}.
The state at time $t$ is now $\xi_t \dvtx C_N (\omega) \longrightarrow\{
1, 2 \}$, that is, each vertex is occupied by a symbiont of one
of two types.
Letting for $i = 1, 2$
\begin{eqnarray*}
f_i (x) & = & \card\{y \dvtx\pi(y) = \pi(x) \mbox{ and } \xi_t
(y) = i \} / N, \\
g_i (x) & = & \card\{y \dvtx\pi(y) \sim\pi(x) \mbox{ and } \xi_t
(y) = i \} / (N \deg\pi(x))
\end{eqnarray*}
denote the fraction of type $i$ symbionts at site $\pi(x)$ and its
neighborhood, respectively, the evolution is described
by the Markov generator $L_2$ defined on the set of the cylinder
functions by
\begin{eqnarray*}
L_2 f (\xi) & = &
\sum_{x \in C_N (\omega)} \frac{\alpha_1 f_1 (x)
+ \beta_1 g_1 (x)}{\alpha_1 f_1 (x) +
\alpha_2 f_2 (x) + \beta_1 g_1 (x) + \beta_2 g_2 (x)} [f (\xi_{x, 1}) -
f (\xi)]
\\ &&{}
+ \sum_{x \in C_N (\omega)} \frac{\alpha_2 f_2 (x)
+ \beta_2 g_2 (x)}{\alpha_1 f_1 (x) +
\alpha_2 f_2 (x) + \beta_1 g_1 (x) + \beta_2 g_2 (x)} [f (\xi_{x, 2}) -
f (\xi)],
\end{eqnarray*}
where $\xi_{x, i}$ is the configuration obtained from $\xi$ by
assigning the value $i$ to vertex~$x$.
The transition rates indicate that, regardless of its type, each
symbiont dies at rate~1 and gets instantaneously replaced
by a symbiont whose type is chosen from the nearby symbionts according
to the relative fecundities and transmissibilities
of the two symbiont species.
In the neutral case when the reproduction rates are both equal to say
$\alpha$ and the transmission rates are both equal to
say $\beta$, the local evolution reduces to the following:
the type of each symbiont is updated at rate 1 and the new type is
chosen uniformly at random from the same host with
probability $\alpha/ (\alpha+ \beta)$ or a~nearby host with
probability $\beta/ (\alpha+ \beta)$.
Note that the process $\{\xi_t \}_t$ can again be seen as a mixture of
two well-known processes, namely, the Moran model with
selection, and its spatial analog, the biased voter model \cite
{bramsongriffeath1980,bramsongriffeath1981} on the infinite
percolation cluster. See Figure \ref{fig:SHMV} for simulation
pictures of this voter model.

%
%
%


To state our results for the competition model, we set $\theta\in(0,
1)$ and denote by $\pi_{\theta}$ the product measure
restricted to $C_N (\omega) = C_{\infty} (\omega) \times\K_N$
defined by
\[
\pi_{\theta} \{\xi(x) = 1 \} = \theta\quad\mbox{and}\quad
\pi_{\theta} \{\xi(x) = 2 \} = 1 - \theta\qquad\mbox{for
all } x \in C_N (\omega).
\]
From now on, we assume that $\{\xi_t \}_t$ starts from the product
measure $\pi_{\theta}$ and let $\Rightarrow$ stand for convergence
in distribution.
The process is said to cluster if
\[
\mbox{there exists } a \in(0, 1) \qquad\mbox{such that } \xi_t
\Rightarrow a \delta_1 + (1 - a) \delta_2 \qquad\mbox{as } t \to
\infty.
\]
In particular, we have
\[
\lim_{t \to\infty} P \bigl(\xi_t (x) \neq\xi_t (y)\bigr) = 0\qquad
\mbox{for all } x, y \in C_N (\omega).
\]
The process is said to coexist if in contrast $\xi_t \Rightarrow\nu
_{\theta}$ as $t \to\infty$ for some $\nu_{\theta}$ such that
\[
\nu_{\theta} \{\xi(x) \neq\xi(y) \} \neq0 \qquad\mbox{for
all } x, y \in C_N (\omega), x \neq y.
\]
Type 1 is said to invade type 2 if
\[
P \Bigl(\lim_{t \to\infty} N_t = \infty\bigm| N_0 = 1 \Bigr) > 0 \qquad\mbox{where } N_t =
\card\{x \in C_N (\omega) \dvtx\xi
_t (x) = 1 \},
\]
indicating that, starting with a single symbiont of type 1 in the
infinite percolation cluster, there is a positive probability that
the number of type~1 keeps growing indefinitely.
Finally, type 1 is said to outcompete type 2 whenever we have the
stronger condition $\xi_t \Rightarrow\delta_1$.
\begin{theor}[(Voter interactions)]
\label{competition}
Assume that $p > p_c$.
If $\alpha_1 = \alpha_2$ and $\beta_1 = \beta_2$, the component
community clusters in two dimensions, whereas
coexistence occurs in higher dimensions.
\end{theor}

The analysis of the neutral competition model relies on duality techniques.
We show, in the neutral case, that the process is dual to a certain
system of coalescing random walks evolving on the random graph
induced by the infinite percolation cluster.
The long-term behavior of the process is related to the so-called
finite/infinite collision property of the graph, which is studied
in details in two dimensions and higher dimensions separately in
Section \ref{sec:competition}.
\begin{theor}[(Voter interactions with selection)]
\label{selection}
Assume that $p = 1$.
\begin{enumerate}
\item If $\alpha_1 \geq\alpha_2$ and $\beta_1 > \beta_2$, then
type 1 invades type 2.
%
\item If $\alpha_2 = 0$ and $\beta_1 > \beta_2$, then type 1
outcompetes type 2.
\end{enumerate}
\end{theor}

Note that, in contrast with the neutral case, when one symbiont type
has a selective advantage, the transition rates of the
competition model are no longer linear functions of the local
frequencies of symbiont types.
This leads a priori to a failure of duality techniques.
The first part of the theorem is established by comparing the number of
type 1 symbionts with a gambler's ruin model whereas
the second part relies on the analysis of a certain system of branching
coalescing random walks which is dual to a particle
system related to the model with selection when $\alpha_2 = 0$.
Before going into the proofs, we note that, though the condition
$\alpha_2 = 0$ may appear biologically unrealistic, the
second part obviously holds when $N = 1$ since in this case the value
of the reproduction parameters becomes unimportant.
In particular, the result in the second part is relevant for species in
which only one symbiont individual can associate with
a host individual.
Finally, we point out that the techniques to prove Theorem \ref
{selection} also apply to the case when $p > p_c$
but lead to conditions on the parameters which are far from being
optimal and to very tedious calculations that only make
the key ideas unclear.
Therefore, for simplicity, we focus on the case $p = 1$ only.

\section{\texorpdfstring{Proof of Theorem \protect\ref{invasion}}{Proof of Theorem 1}}
\label{sec:invasion}

This section is devoted to the analysis of the invasion model, and more
precisely to the proof\vadjust{\goodbreak} of the second part
of Theorem \ref{invasion}.
The key idea is to show that the branching random walk restricted to a
large square persist an arbitrary long
time provided $\alpha+ \beta> 1$ and $N$ is large.
The combination of our estimates with a block construction implies
survival of the metapopulation restricted to an infinite
self-avoiding path of large squares fully occupied by hosts.
It is also proved that such a path exists whenever the parameter $p$ is
close enough to 1.

\begin{figure}[b]

\includegraphics{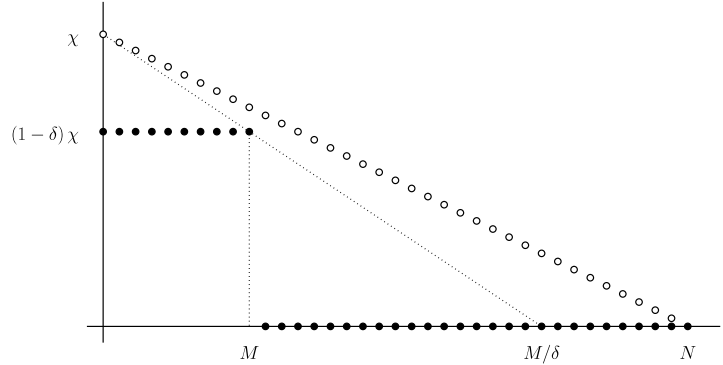}

\caption{Reproduction rates for $\{\zeta_t \}_t$ and $\{\bar\eta_t
\}_t$ ($\bullet$ and $\circ$, resp.).}
\label{fig:domination}
\end{figure}

The first step is to prove branching random walk estimates in order to
establish the result in any dimension when $p = 1$,
that is, the infinite percolation cluster consists of the entire lattice.
To begin with, we observe that, for all $M > 0$ and $\delta\in(0,
1)$, the process $\{\bar\eta_t \}_t$ dominates, for $N$
sufficiently large, the process $\{\zeta_t \}_t$ whose dynamics are
described by the Markov generator
\begin{eqnarray*}
D_1 f (\zeta) &=& \sum_{X \in C_{\infty} (\omega)}
\zeta(X) [f (\zeta_{X-}) - f (\zeta)] \\
&&{}+ \sum_{X \in C_{\infty} (\omega
)} (1 - \delta) \ind\{\zeta(X) \leq M \}
\biggl(\alpha\zeta(X) + \frac{\beta}{2d} \sum_{X \sim Y} \zeta(Y)
\biggr)\\
&&\hspace*{47.6pt}{}\times [f (\zeta
_{X+}) - f (\zeta)],
\end{eqnarray*}
where the configurations $\zeta_{X-}$ and $\zeta_{X+}$ are obtained
from $\zeta$ by, respectively, removing and adding a symbiont
at site $X$.
Indeed, it suffices that $N \geq M / \delta$ since in that case
\[
1 - j N^{-1} \geq1 - \delta j M^{-1} \geq(1 - \delta) \ind
\{j \in[0, M] \} \qquad\mbox{for all } j = 0, 1, \ldots, N.
\]
To see this, we observe that the process $\{\zeta_t\}_{t}$ is a
truncated branching random walk that allows at most $M + 1$ particles
per site at the same time.
See Figure \ref{fig:domination} where we compare the reproduction
rates to site $Y$ of a particle living at site $X$ for the
processes $\{\bar\eta_t\}_t$ and $\{\zeta_t\}_t$: on the $x$-axis we
have the number of particles at $Y$, and the parameter $\chi$ is
equal to $\alpha$ if $X = Y$ and~$\beta$ if $X \sim Y$.

%

Let $\delta> 0$ such that $(\alpha+ \beta)(1 - \delta) > 1$.
We will prove, following \cite{bertacchizucca2009}, that for all $M$
sufficiently large, the truncated branching
random walk $\{\zeta_t \}_t$ survives with positive probability, by
looking at the branching random walk $\{\bar\zeta_t \}_t$
whose dynamics are described by
\begin{eqnarray*}
\bar D_1 f (\bar\zeta) &=& \sum_{X \in C_{\infty}
(\omega)} \bar\zeta(X) [f (\bar\zeta_{X-}) - f (\bar\zeta)] \\
&&{}+ \sum_{X \in C_{\infty} (\omega)} \biggl(\bar\alpha
\bar\zeta(X) + \bar\beta\sum_{X \sim Y} \bar\zeta(Y) \biggr) \\
&&\hspace*{48.5pt}{}\times[f (\bar\zeta
_{X+}) - f (\bar\zeta)]
\end{eqnarray*}
starting with one particle at the origin, where $\bar\alpha= (1 -
\delta) \alpha$ and $\bar\beta= (1 - \delta) \beta/ 2d$.
\begin{lemma}
\label{BRW}
For $X \!\sim\!0$, we have $\E[\bar\zeta_n (X) \mid \bar\zeta_0
(0) \!=\! 1] \!>\! 1$ for $n$ large enough.
\end{lemma}
\begin{pf}
We observe that $\E(\bar\zeta_t (X))$ satisfies the differential
equation\break (see~\cite{belhadjibertacchizucca2009}, Section
4)
\[
\frac{d}{dt} \E(\bar\zeta_t (X)) = - \E(\bar\zeta_t
(X)) + \bar\alpha\E(\bar\zeta_t (X))
+ \bar\beta\sum_{X \sim Y} \E(\bar\zeta_t (Y)),
\]
whose solution is
%
%
\begin{equation}
\label{eq:BRW1}
\E(\bar\zeta_t (X)) = \sum_{n = 0}^{\infty} \sum_{k = 0}^{n - 1} \mu
^{(n, k)} (0, X) \frac{\bar\alpha^k \bar\beta^{n - k} t^n}{n!} e^{-t},
\end{equation}
where $\mu^{(n, k)} (0, X)$ is the number of paths from site 0 to site
$X$ of length $n$ with $k$ loops.
To estimate the right-hand side of (\ref{eq:BRW1}), we let $\{
U_k \}_k$ be the discrete-time random walk with
\[
P (U_{k + 1} = Z \mid U_k = Y) = \cases{
\bar\alpha(\bar\alpha+ 2d \bar\beta)^{-1}, &\quad for $Y =
Z$,\cr
\bar\beta(\bar\alpha+ 2d \bar\beta)^{-1}, &\quad for $Y \sim Z$,}
\]
and observe that, for any site $X \sim0$,
\begin{eqnarray*}
\sum_{k = 0}^{n - 1} \mu^{(n, k)} (0, X) \frac{\bar\alpha^k
\bar\beta^{n - k}}{(\bar\alpha+ 2d \bar\beta)^n} &=& P (U_n = X \mid U_0 =
0) \\[-1pt]
&\geq& C_1 n^{-d/2}
\end{eqnarray*}
for a suitable $C_1 = C_1 (\bar\alpha, \bar\beta) > 0$.
We refer to \cite{woess2000}, Corollary 13.11, for the asymptotic
estimates of the $n$-step probabilities.
In particular,\vadjust{\goodbreak} for $X \sim0$ and $t = n$, we obtain
\begin{eqnarray*}
\E(\bar\zeta_n (X)) &\geq& \sum_{k = 0}^{n - 1} \mu^{(n, k)} (0, X)
\frac
{\bar\alpha^k \bar\beta^{n - k} n^n}{n!} e^{-n}
\\[-1pt]
&=& \sum_{k = 0}^{n - 1} \mu^{(n, k)} (0, X) \frac
{\bar\alpha^k \bar\beta^{n - k}}{(\bar\alpha+ 2d \bar\beta
)^n} \frac{n^n (\bar\alpha+ 2d \bar\beta)^n}{n!} e^{-n}
\\[-1pt]
&\stackrel{n \to\infty}{\sim}& \frac{(\bar\alpha+ 2d \bar\beta)^n}{\sqrt
{2 \pi n}}
\sum_{k = 0}^{n - 1} \mu^{(n, k)} (0, X) \frac{\bar\alpha^k \bar\beta
^{n - k}}{(\bar\alpha+ 2d \bar
\beta)^n} \\[-1pt]
&\geq&\frac{(\bar\alpha+ 2d \bar\beta)^n}{\sqrt{2 \pi n}} C_2 n^{-d/2}
\end{eqnarray*}
for a suitable $C_2 > 0$.
Finally, since $(\bar\alpha+2d \bar\beta)=(\alpha+ \beta) (1 -
\delta) > 1$, we deduce that
\[
\E(\bar\zeta_n (X)) \geq\frac{(1 - \delta)^n (\alpha+
\beta)^n}{\sqrt{2 \pi n}} C_2 n^{-d/2} > 1
\]
provided $n$ is sufficiently large.\vadjust{\eject}
\end{pf}

Following the ideas of Lemma 5.3, Remark 5.2 and Theorem 5.1
in \cite{bertacchizucca2009}, and using Lemma \ref{BRW}
above in place of \cite{bertacchizucca2009}, Lemma 5.2, one proves
that $\{\zeta_t \}_t$ survives when $M$ is sufficiently
large, and so does, by stochastic domination, the metapopulation when
$N$ is large and the density $p = 1$.
These ideas are developed in more details in the following lemma.
\begin{lemma}
\label{survival0}
If $(\alpha+ \beta) (1 - \delta) > 1$ then the process $\{\zeta_t\}
_t$ survives when $M$ is sufficiently large.
\end{lemma}
\begin{pf}
By additivity of $\{\bar\zeta_t \}_t$, if $X \sim0$ then the central
limit theorem implies that
\[
\lim_{K \to\infty} \biggl[\!
P \bigl(\bar\zeta_n(X) \geq K \mid \bar\zeta_0(0) = K\bigr) - 1 + \Phi\biggl(\!\frac{K - \E
(\bar\zeta_n(X) \mid \bar\zeta_0 (0) =
1) K}{\sqrt{{\mathrm{Var}}(\bar\zeta_n (X) \mid \bar\zeta_0
(0) = 1)} \sqrt K} \!\biggr)\!
\biggr] = 0,
\]
where the function $\Phi$ is the cumulative distribution function of
the standard normal.
Since $n$ is fixed, it follows that, for all $\ep> 0$,
\[
P \bigl(\bar\zeta_n (X) \geq K \mbox{ for all } X \sim0 \mid \bar
\zeta_0 (0) = K\bigr) > 1 - \ep
\]
for $K$ sufficiently large.
Let $\{N_t \}_t$ be the branching process with birth rate $\bar\alpha
+ 2d \bar\beta$ and death rate zero, which represents
the total number of particles born up to time~$t$.
By the same argument as before, there exists $C_3 > 1$ such that
\[
P (N_n \leq C_3 K \mid N_0 = K) \geq1 - \ep\qquad\mbox{for
all $K$ sufficiently large}.
\]
Since, if $M \geq C_3 K$ then $\{\bar\zeta_t\}_t$ and $\{\zeta_t\}
_t$ coincide (up to time $n$) on $\{N_n \leq C_3 K\}$,
we have
\[
P \bigl(\zeta_n (X) \geq K \mid \zeta_0 (0) = K\bigr) > 1 - 2 \ep.
\]
In order to get
%
%
\begin{equation}
\label{eq:BRW2}
P \bigl(\zeta_n (X) \geq K \mbox{ for all } X \sim0 \mid \zeta_0
(0) = K\bigr) > 1 - 2 \ep
\end{equation}
we need to ensure that from time 0 to time $n$, in no site the process
$\{\bar\zeta_t\}_t$ on $\{N_n \leq C_3 K \}$
ever exceeds $M$ particles.
By geometric arguments (see~\cite{bertacchizucca2009}, Step 3, for
further details), one proves that it suffices to take
$M \ge2 H_0 C_3 K = C_4 K$ where $H_0$ is the number of paths of
length $n$ in $\Z^d$ crossing a fixed vertex.
To complete the proof, we couple the process $\{\zeta_t \}_t$ with a
supercritical 1-dependent oriented site percolation process
on $\Z\times\Z_+$ in a way such that the existence of an infinite
cluster implies survival for~$\{\zeta_t \}_t$ relying on
the standard rescaling technique introduced in \cite
{bramsondurrett1988}. Let
\[
\mathcal G = \{(z, m) \in\Z\times\Z_+ \dvtx z + m \mbox{ is even}
\},
\]
and declare site $(z, m) \in\mathcal G$ to be good if the host at site
$z e_1$ is infected by at least $K$ symbionts at time
$m \times n$, where $e_1$ denotes the first unit vector.
Also, let
\[
\mathcal G_m = \{z \in\Z\dvtx(z, m) \mbox{ is a good site} \}
\]
denote the set of good sites at level $m$.
Then, inequality (\ref{eq:BRW2}) above implies that $\mathcal G_m$
dominates stochastically the set $W_m$ of wet sites at level
$m$ in a 1-de\-pendent oriented site percolation process on the lattice
$\mathcal G$ with parameter $1 - \ep$ and with initial
condition $W_0 \subset\mathcal G_0$ (see Durrett \cite{durrett1984}
for a complete description of oriented percolation).
The result then follows by choosing \mbox{$\ep>0$} sufficiently small to
make the oriented percolation process supercritical.
\end{pf}

Since, in the proof of Lemma \ref{BRW}, we consider only the
particles of genera\-tion~$n$, equation (\ref{eq:BRW2}) holds
if, instead of the process $\{\zeta_t \}_t$, we deal with the process
$\{\zeta^n_t \}_t$ obtained by deleting all the particles of
generation $n' > n$.
In addition, the process $\{\bar\eta_t^n \}_t$, obtained from the
metapopulation model by assuming that symbionts sent
outside $[n, n]^d$ are killed, clearly dominates~$\{\zeta^n_t \}_t$.
\begin{lemma}
\label{survival}
Fix $n$ so that Lemma \ref{BRW} holds.
Then, for all $\ep> 0$,
\[
P \bigl(\bar\eta_n^n (X) \geq\sqrt N \mbox{ for all } X \sim0 \mid \bar\eta_0^n
(0) \geq\sqrt N\bigr) \geq1 - \ep
\]
for all $N$ sufficiently large.
\end{lemma}
\begin{pf}
This follows from (\ref{eq:BRW2}) (using $\{\bar\eta_t^n \}
_t$ instead of $\{\zeta_t^n \}_t$) choosing $K=\sqrt{N}$,
from stochastic domination when $C_4 \sqrt N / \delta< N$, and from
the monotonicity of $\{\bar\eta_t^n \}_t$.
\end{pf}

\begin{figure}

\includegraphics{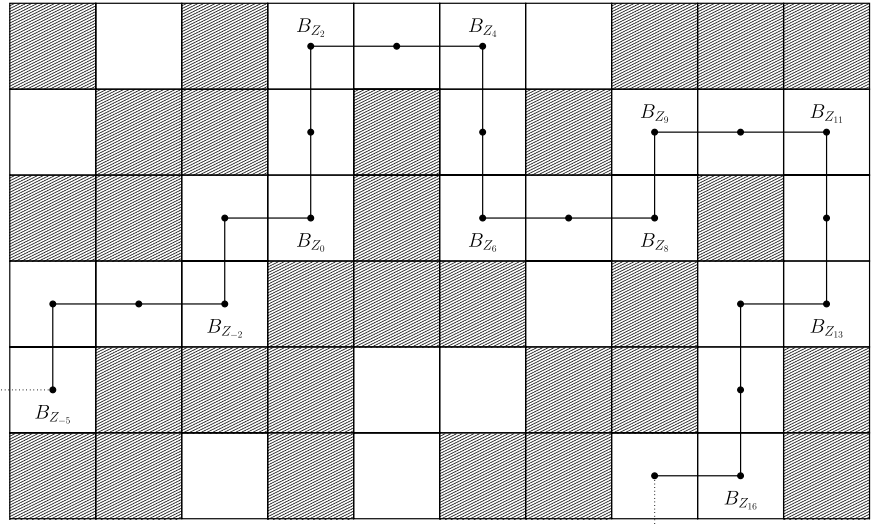}

\caption{Picture of the self-avoiding path $\Gamma$.}
\label{fig:path}
\end{figure}

To deduce the second part of Theorem \ref{invasion} from the
previous lemma, we use another block construction in order
to compare the evolution of the metapopulation along an infinite
self-avoiding path with oriented percolation.
To apply successfully Lemma \ref{survival}, any site within distance
$n$ of this infinite self-avoiding path must be open.
The existence of such a path follows by choosing the parameter $p$
close enough to 1.
First, we fix $n$ so that Lemma \ref{BRW} holds (recall that $n$ only
depends on the reproduction rate $\alpha$, the transmission
rate $\beta$, and the spatial dimension $d$).
Then, we fix the parameter $\ep> 0$ such that $1 - \ep$ is greater
than the critical value of $n$-dependent oriented percolation.
We prove the result when the density $p$ of hosts satisfies
\[
p > \exp\bigl((2n + 1)^{-d} \log p_c\bigr),
\]
where $p_c$ is the critical value of site percolation in $d$ dimensions.
We tile the $d$-dimensional regular lattice with cubes of edge length
$2n + 1$ by setting
\[
B_0 = [- n, n]^d \quad\mbox{and}\quad B_Z = (2n + 1) Z +
B_0 \qquad\mbox{for all } Z \in\Z^d.
\]
Given a realization $\omega$ of the site percolation process with
parameter $p$, we call a cube $B_Z$ open if all the sites $X \in B_Z$
are occupied by a host, and closed otherwise.
Our choice of $p$ implies
\[
P (B_Z \mbox{ is open}) = p^{(2n + 1)^d} > p_c \qquad\mbox
{for all } Z \in\Z^d.
\]
In particular, there exists almost surely an infinite self-avoiding
path of open cubes, that is, there exists a self-avoiding
path $\{Z_i \dvtx i \in\Z\} \subset\Z^d$ such that cube $B_{Z_i}$ is
open for all $i$.
From this path of open cubes, we construct an infinite self-avoiding
path of open sites $\Gamma= \{\Gamma_z \dvtx z \in\Z\}$ by including
all the sites belonging to the straight lines connecting the centers of
adjacent cubes, as shown in Figure \ref{fig:path} where gray
%
%
squares refer to closed cubes, and white squares to open cubes.
By construction:
\begin{enumerate}
\item For all $z \in\Z$ and all $X \in\Gamma_z + [- n, n]^d$, we
have $X \in C_{\infty} (\omega)$.
\item For all $z \in\Z$, we have $\Gamma_z \sim\Gamma_{z + 1}$.
\end{enumerate}
Site $(z, m) \in\mathcal G$ is now said to be good whenever the host
at $\Gamma_z$ is infected by at least $\sqrt N$ symbionts at time
$m \times n$.
As previously, we let $\mathcal G_m$ denote the set of good sites at
level $m$.
Then Lemma \ref{survival} and the fact that the evolution rules of the
process are homogeneous in time imply that
%
%
\begin{equation}
\label{eq:BRW3}
P \bigl((z - 1, m + 1) \mbox{ and } (z + 1, m + 1) \mbox{ are good} \mid (z, m)
\mbox{ is good}\bigr) \geq1 - \ep
\end{equation}
for sufficiently large $N$.
Denoting again by $W_m$ the set of wet sites at level $m$ in an
$n$-dependent oriented site percolation process with parameter $1 - \ep$
the inequality (\ref{eq:BRW3}) implies that the processes can be
constructed on the same probability space in such a way that
\[
P (W_m \subset\mathcal G_m \mbox{ for all } m \geq0 \mid W_0
\subset\mathcal G_0) = 1.
\]
Since $1 - \ep$ is greater than the critical value of oriented
percolation, this implies as previously that the metapopulation survives,
which completes the proof of Theorem \ref{invasion}.

\section{\texorpdfstring{Proof of Theorem \protect\ref{competition}}{Proof of Theorem 2}}
\label{sec:competition}

This section is devoted to the analysis of the competition model under
neutrality.
The process can be constructed graphically relying on an idea of Harris
\cite{harris1972} from a collection of independent Poisson
processes.
In the neutral case, because the transition rates are linear functions
of the local frequencies, the graphical representation
induces a natural duality relationship between the spatial model and a
system of coalescing random walks on $C_N (\omega)$, and
Theorem \ref{competition} follows from certain collision properties of
symmetric random walks on the infinite percolation
cluster.


\subsection*{Duality with coalescing random walks}

To define the dual process of the competition model under neutrality,
we first construct the process graphically
from collections of independent Poisson processes using an idea of
Harris \cite{harris1972}.
Each vertex $x \in C_N (\omega)$ is equipped with a Poisson process
with parameter 1.
Poisson processes attached to different vertices are independent.
At the arrival times of the process at $x$, we toss a coin with success
probability $\alpha/ (\alpha+ \beta)$
where $\alpha$ is the common reproduction parameter of both symbiont
types and $\beta$ the common transmission parameter.
If there is a success, we choose a vertex uniformly at random from the
host at site~$\pi(x)$ and draw an arrow
from this vertex to vertex $x$.
If there is a failure, we choose a vertex uniformly at random from one
of the hosts adjacent to site $\pi(x)$ and draw
an arrow from this vertex to vertex $x$.
In view of the geometry of the graph and the number of vertices per
host, this is equivalent to saying that:
\begin{itemize}
\item[--] For any pair of vertices $x, y \in C_N (\omega)$ with $x
\ver y$, we draw an arrow from $y$ to $x$ at the
arrival times of an independent Poisson process with parameter $\alpha
/ (N (\alpha+ \beta))$.
\item[--] For any pair of vertices $x, y \in C_N (\omega)$ with $x
\hor y$, we draw an arrow from $y$ to $x$ at the
arrival times of an independent Poisson process with parameter $\beta/
(N \deg\pi(x) (\alpha+ \beta))$.
\end{itemize}
In any case, an arrow from vertex $y$ to vertex $x$ indicates that the
symbiont at $x$ dies and gets instantaneously
replaced by a symbiont of the same species as the symbiont at $y$.

\begin{figure}

\includegraphics{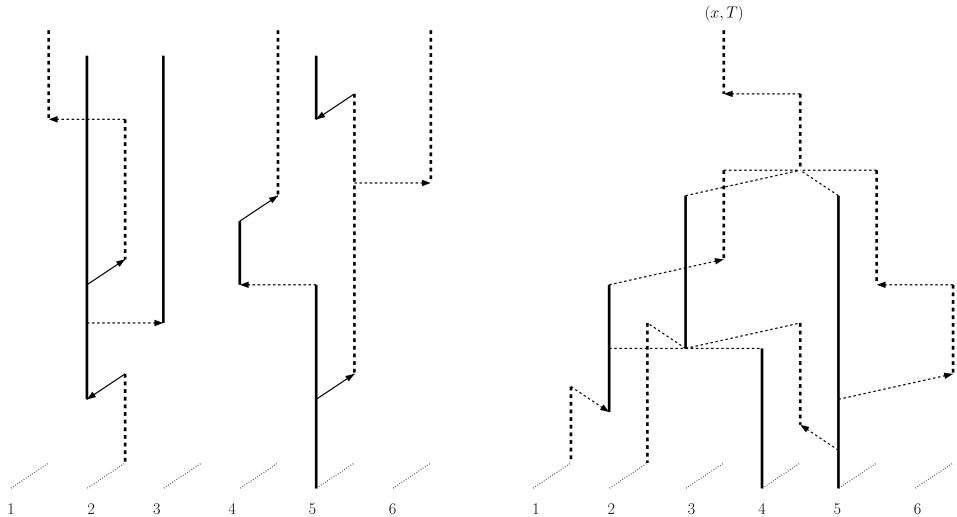}%
\vspace*{-3pt}
\caption{Dual process in the neutral case and branching random walk
$\hat\zeta_s (x, T)$.
In both pictures, $N = 2$ and for simplicity we have set $p = 1$ and $d
= 1$.
Time goes up, and arrows within the same host are drawn in continuous
line, while arrows connecting two adjacent hosts are drawn in dashed lines.}
\label{fig:voter}
\vspace*{-5pt}
\end{figure}

To define the dual process, we say that there is a path from $(y, T -
s)$ to $(x, T)$, which corresponds to a
dual path from $(x, T)$ to $(y, T - s)$, if there are sequences of
times and vertices
\[
s_0 = T - s < s_1 < \cdots< s_{n + 1} = T
\quad\mbox{and}\quad
x_0 = y, x_1, \ldots, x_n = x
\]
such that the following two conditions hold:
\begin{enumerate}
\item for $i = 1, 2, \ldots, n$, there is an arrow from $x_{i - 1}$ to
$x_i$ at time $s_i$ and
\item for $i = 0, 1, \ldots, n$, there is no arrow that points at the
segments $\{x_i \} \times(s_i, s_{i + 1})$.
\end{enumerate}
The dual process starting at $(x, T)$ is the process defined by
\[
\hat\xi_s (x, T) = \{y \in C_N (\omega) \dvtx\mbox{there is a dual
path from $(x, T)$ to $(y, T - s)$} \}.
\]
The dual process starting from a finite set of vertices $B \subset C_N
(\omega)$ can be defined as well.
In this case, the dual process starting at $(B, T)$ is the set-valued
process defined by
\begin{eqnarray*}
\hat\xi_s (B, T) & = &
\{y \in C_N (\omega) \dvtx\mbox{there is a dual path} \\
&&\hspace*{5.1pt}
\mbox{from $(x, T)$ to $(y, T - s)$ for some $x \in B$}
\} \\ & = &
\{y \in C_N (\omega) \dvtx y \in\hat\xi_s (x, T) \mbox{ for some } x
\in B \}.
\end{eqnarray*}
The dual process is naturally defined only for $0 \leq s \leq T$.
However, it is convenient to assume that the Poisson processes in the
graphical representation are defined for
negative times so that the dual process can be defined for all $s \geq0$.
Note that, in view of the graphical representation of the competition
model in the neutral case, the dual process
starting at $(x, T)$ performs a continuous-time random walk on the
random graph $C_N (\omega)$ that makes transitions
\[
y \to\cases{
z, &\quad for $z \ver y$ at rate $\alpha/ \bigl(N (\alpha+
\beta)\bigr)$,\cr
z, &\quad for $z \hor y$ at rate $\beta/ \bigl(N \deg\pi(y) (\alpha+ \beta
)\bigr)$.}
\]
The dual process starting from a finite set $B \subset C_N (\omega)$
consists of a system of $\card(B)$ such random
walks, one random walk starting from each vertex in the set $B$.
Any two of these random walks evolve independently of each other until
they intersect when they coalesce.
This induces a duality relationship between the model and coalescing
random walks.
We refer the reader to the left-hand side of Figure \ref{fig:voter}
for an example of realization of the dual process
in the neutral case.

%

The reason for introducing the dual process is that it allows us to
deduce the configuration of the system
at the current time based on the configuration at earlier times, but
also how vertices at the current time are correlated,
by keeping track of the ancestry of each symbiont.
In particular, the long-term behavior of the competition model
(clustering versus coexistence) can be expressed in terms
of collision properties of random walks on the infinite percolation
cluster through the duality relationship between the
model and coalescing random walks.
We now explain this connection in details, starting with some key definitions.

Let $G = (V, E)$ be an infinite connected graph.
We call simple symmetric random walk on this graph the continuous-time
Markov process $\{X_t \}_t$ with state space $V$
that jumps from $u$ to $v$ at rate one if and only if $(u, v) \in E$.
Note that the embedded Markov chain associated to this Markov process
is the discrete-time random walk $\{\X_n \}_n$
with transition probabilities
\[
P (\X_{n + 1} = v \mid \X_n = u) = \frac{1}{\deg(u)}
\quad\mbox{if and only if}\quad (u, v) \in E.
\]
Since the graph $G$ is connected, the process $\{X_t \}_t$ is
irreducible so either all the vertices of the graph are
recurrent, in which case the graph is said to be recurrent, or all the
vertices are transient, in which case the graph
is said to be transient.
Let $X_t$ and $Y_t$ be two independent random walks on the graph~$G$.
Using again the fact that the graph is connected and the Kolmogorov
zero-one law, the probability that the two random
walks intersect infinitely often, namely
\[
P (\mbox{for all $t$ there exists $s$ such that $X_{t + s} = Y_{t +
s}$})
\]
is either equal to 0 or 1 regardless of the initial positions of the
random walks.
The graph $G$ is said to have the infinite collision property if the
previous probability is equal to 1, and it is
said to have the finite collision\vadjust{\eject} property if the previous probability
is equal to 0.
Such properties for the infinite percolation cluster $C_{\infty}
(\omega)$ translate through the duality relationship
with coalescing random walks into coexistence/clustering of the
competition model, as shown in the next lemma.
\begin{lemma}
\label{duality}
We have the following alternative:
\begin{enumerate}
\item$C_{\infty} (\omega)$ has the infinite collision property and
then the process clusters, or
\item$C_{\infty} (\omega)$ has the finite collision property and
then coexistence occurs.
\end{enumerate}
\end{lemma}
\begin{pf}
Let $B \subset C_N (\omega)$ be finite and let $\Theta^i_t = \{x \in
C_N (\omega) \dvtx\xi_t (x) = i \}$ denote the set
of vertices occupied by a type $i$ symbiont at time $t$. By duality,
%
%
\begin{equation}
\label{eq:dual1}
P (\Theta^1_t \cap B = \varnothing) = \E(1 - \theta)
^{|\hat\xi_t (B, t)|},
\end{equation}
where $\theta$ is the initial density of type 1.
Since the number of particles in $\hat\xi_t (B, t)$ is a
nonincreasing function of $t$ and has a limit, the bounded
convergence theorem implies that the probability on the left-hand side
of (\ref{eq:dual1}) also has a limit
as $t \to\infty$.
It follows that the process converges to a stationary distribution.
To understand how different vertices are correlated under this
stationary distribution, we take two vertices
$x, y \in C_N (\omega)$, $x \neq y$, and consider the projections
\[
X_s = \pi(\hat\xi_s (x, T)) \quad\mbox{and}\quad Y_s = \pi(\hat\xi_s (y, T)).
\]
Let $\tau$ be the hitting time of the dual processes, that is,
\[
\tau= \inf\{s > 0 \dvtx\hat\xi_s (x, T) = \hat\xi_s (y, T) \}.
\]
Note that the processes $X_s$ and $Y_s$ evolve individually according
to continuous-time random walks run at rate
$q := \beta/ (\alpha+ \beta)$ on the infinite percolation cluster
$C_{\infty} (\omega)$.
They evolve independently of each other until time $\tau$ when they coalesce.
We set $t_0 = 0$ and define inductively
\begin{eqnarray*}
s_i & = & \inf\{s > t_{i - 1} \dvtx X_s = Y_s \}, \\
t_i & = & \inf\{s > s_i \dvtx X_s \neq Y_s \} = \inf\{s > s_i \dvtx
X_s \neq X_{s_i} \mbox{ or } Y_s \neq Y_{s_i} \}
\end{eqnarray*}
for $i \geq1$.
Note that, if the dual processes coalesce at time $s$ then
\[
t_i = s_{i + 1} = \infty\qquad\mbox{for all } i \geq
\max\{j \dvtx s_j < s \}.
\]
Also, let $M_i$ denote the total number of jumps during the interval of
time $(s_i, t_i)$ of either of the dual processes
starting at vertex $x$ or vertex $y$.
Wri\-ting~$P_i$ for the conditional probability given the event that $s_i
< \infty$ and using that each dual process jumps
to one of the adjacent hosts at rate $q$ and within\vadjust{\eject} each host at rate
$1 - q$, we obtain the following probability:
%
%
\begin{eqnarray}
\label{eq:dual2}
P_i (\tau> t_i) &=& \sum_{j = 0}^{\infty} P_i (\tau> t_i \mid M_i =
j) P_i (M_i = j) \nonumber\\[-3pt]
&=&
\sum_{j = 0}^{\infty} \biggl(1 - \frac{1}{N}
\biggr)^j q (1 - q)^j = \frac{q}{1 - (1 - q) (1 - 1 / N)} \\[-3pt]
&=& \frac{N \beta
}{\alpha+ N \beta}.
\nonumber
\end{eqnarray}
Let $J = J (x, y, \omega) = \sup\{j \dvtx s_j < \infty\}$, and note
that, on the event that the dual processes starting at
$x$ and $y$ do not coalesce, $J$ is equal in distribution to the number
of intersections of two independent random walks on
the infinite percolation cluster starting at $\pi(x)$ and $\pi(y)$.
In particular, if~$C_{\infty} (\omega)$ has the infinite collision
property and if $I_t$ denotes the number of intersections
up to time $t$ of two independent random walks starting at $\pi(x)$
and~$\pi(y)$ then (\ref{eq:dual2}) implies
\begin{eqnarray*}
\lim_{T \to\infty} P \bigl(\xi_T (x) \neq\xi_T (y)\bigr)
& \leq&
\lim_{T \to\infty} P \bigl(\hat\xi_T (x, T) \neq\hat
\xi_T (y, T)\bigr)
\\[-3pt] & = &
\lim_{T \to\infty} \sum_{j = 0}^{\infty} \prod
_{i = 0}^j P (I_T = j) P_i (\tau> t_i)
\\[-3pt] & = &
\lim_{T \to\infty} \sum_{j = 0}^{\infty} \biggl(\frac{N \beta}{\alpha+ N
\beta} \biggr)^j P (I_T = j) = 0
\end{eqnarray*}
by the bounded convergence theorem since $P (I_T = j) \to0$ as $T
\to\infty$ for all $j \in\N$.
This shows that the process clusters.
Alternatively, if the infinite percolation cluster has the finite
collision property, then $J$ is almost surely finite so
(\ref{eq:dual2}) implies that
\begin{eqnarray*}
&&
\lim_{T \to\infty} P \bigl(\xi_T (x) \neq\xi_T (y)\bigr) \\[-3pt]
&&\qquad= \theta(1 - \theta)
\lim_{T \to\infty} P \bigl(\hat
\xi_T (x, T) \neq\hat\xi_T (y, T)\bigr)
\\[-3pt]
&&\qquad= \theta(1 - \theta) \sum_{j = 0}^{\infty} \lim_{T \to\infty} P \bigl(\hat\xi
_T (x, T) \neq\hat\xi_T (y, T) \mid J = j\bigr) P (J = j)
\\[-3pt]
&&\qquad= \theta(1 - \theta) \sum_{j = 0}^{\infty} \biggl(\frac{N \beta}{\alpha+ N
\beta} \biggr)^j P (J = j) \geq c > 0,
\end{eqnarray*}
which shows that coexistence occurs.
This completes the proof.
\end{pf}

Note that the previous lemma easily extends to any connected graph in
which the degree of each vertex is uniformly bounded.
That is, given such a graph $G = (V, E)$, the competition model can be
naturally defined on the\vadjust{\eject} graph $G_N$ with vertex set
$V \times\K_N$ constructed from $G$ in the same way as the graph $C_N
(\omega)$ is constructed from the infinite
percolation cluster.
Then, the proof of the previous lemma implies that, for all $N$ finite,
the resulting process clusters when $G$ has the
infinite collision property, but coexists when~$G$ has the finite
collision property.

\vspace*{-3pt}
\subsection*{Transience of the percolation cluster and coexistence}

Motivated by Lem\-ma~\ref{duality}, we now prove that the infinite
percolation cluster $C_{\infty} (\omega)$ in dimensions
$d \geq3$ has the finite collision property.
This follows from the fact that the infinite cluster is transient, a
result due to Grimmett, Kesten and
Zhang~\cite{grimmettkestenzhang1993}, and that the degree of each
vertex is uniformly bounded.
We also answer the same questions for the infinite percolation clusters
in 2 dimensions since the proofs are similar, even though
this result will not be used to establish clustering of the process.\vspace*{-3pt}
\begin{lemma}
\label{walk}
$\!\!$The cluster $C_{\infty} (\omega)$ is recurrent in $d \!=\! 2$ and
transient in $d \!\geq\!3$.\vspace*{-3pt}
\end{lemma}
\begin{pf}
Both statements follow from the fact that any subgraph of a~recurrent
graph is recurrent, and equivalently any supergraph of
a transient graph is transient.
This is proved in details in Doyle and Snell \cite{doylesnell1984}
based on the analogy between random walks and electrical networks
so we only give an outline of their proof.
The idea is to turn the graph under consideration into an electrical
network in which each edge has unit resistance.
Then, simple random walks on this graph are recurrent if and only if
the effective resistance of the resulting electrical
network between a given point and the points at infinity is infinite,
as shown in \cite{doylesnell1984}.
In other respects, Rayleigh's monotonicity law states that if the
resistances of a circuit are increased, respectively, decreased,
then the effective resistance between any two points can only increase,
respectively, decrease.
In particular, removing edges induces an increase of the effective
resistance between any two points, therefore any subgraph of
a recurrent graph is recurrent.
Again, we refer to~\cite{doylesnell1984} for the details.

With the previous result in hands, the first statement follows directly
from the fact that the infinite percolation cluster
in two dimensions is a~subgraph of the two-dimensional lattice, which
is recurrent by Polya's theorem.
Transience of the infinite percolation cluster in higher dimensions has
been established by Grimmett, Kesten
and Zhang \cite{grimmettkestenzhang1993}.
Their proof relies on the construction of a transient tree-like graph
that can be embedded in the infinite percolation cluster.
Note that their result applies to bond percolation but relies on
geometric properties that are known for site percolation as well,
so their proof easily extends to our context.\vspace*{-3pt}
\end{pf}
\begin{lemma}
\label{intersection}
Let $\{X_t \}_t$ and $\{Y_t \}_t$ be two independent random walks run
at rate~1 on $C_{\infty} (\omega)$ both starting at
vertex $A$, and denote by $I (X, Y)$ the number of their intersections. Then,
\[
\E I (X, Y) = \infty\qquad\mbox{in } d = 2 \quad\mbox{and}\quad
\E I (X, Y) < \infty\qquad\mbox{in } d \geq3.
\]
\end{lemma}
\begin{pf}
Since the total rate of jump of both random walks equals 2,
%
%
\begin{eqnarray}
\label{eq:coex1}
\E I (X, Y) &=& 2 \E\biggl(\int_0^{\infty} \ind\{X_t = Y_t \}
\,dt \biggr) = 2 \int_0^{\infty} P (X_t = Y_t) \,dt
\nonumber\\[-2pt]
&=& 2 \int_0^{\infty} \sum_{B \in C_{\infty} (\omega)} P
(X_t = B) P (Y_t = B) \,dt \\[-2pt]
&=& 2 \int_0^{\infty} \sum_{B \in C_{\infty}
(\omega)}
(p_t (A, B))^2 \,dt,\nonumber
\end{eqnarray}
where $p_t (A, B) = P (X_t = B \mid X_0 = A)$.
Now, we observe that the probability that a random walk follows a given
directed path from vertex $A$ to vertex~$B$ is
equal to~1 divided by the product of the degrees of the vertices of
this path excluding the final vertex $B$.
Similarly, the probability that a random walk follows the reverse path
from vertex $B$ to vertex $A$ is 1 divided by the product
of the degrees of the vertices excluding the final vertex $A$, from
which we deduce that
%
%
\begin{equation}
\label{eq:coex2}
(2d)^{-1} p_t (B, A) \leq p_t (A, B) \leq2d p_t (B, A)\qquad
\mbox{for all } B \in C_{\infty} (\omega)
\end{equation}
since $1 \leq\deg(A), \deg(B) \leq2d$.
Therefore, when $C_{\infty} (\omega)$ is recurrent, (\ref
{eq:coex1}) and~(\ref{eq:coex2}) imply that
\begin{eqnarray*}
\E I (X, Y) &\geq& d^{-1} \int_0^{\infty} \sum_{B \in C_{\infty} (\omega
)} p_t (A, B)
p_t (B, A) \,dt \\[-2pt]
&=& d^{-1} \int_0^{\infty} p_{2t} (A, A) \,dt = \infty,
\end{eqnarray*}
whereas when $C_{\infty} (\omega)$ is transient, (\ref
{eq:coex1}) and (\ref{eq:coex2}) imply that
\begin{eqnarray*}
\E I (X, Y) &\leq&4d \int_0^{\infty} \sum_{B \in C_{\infty} (\omega)}
p_t (A, B)
p_t (B, A) \,dt \\[-2pt]
&=& 4d \int_0^{\infty} p_{2t} (A, A) \,dt < \infty.
\end{eqnarray*}
The result then follows from Lemma \ref{walk}.
\end{pf}

Lemma \ref{intersection} indicates that $P (I (X, Y) <
\infty) = 1$ in dimensions $d \geq3$, that is, $C_{\infty} (\omega)$
has the finite collision property, which, together with Lemma~\ref
{duality}, implies that coexistence occurs.
However, that the expected number of intersections is infinite does not
imply that the number of intersections is infinite with
positive probability (with probability 1 by the Kolmogorov zero-one law).
In fact, it is known that recurrent graphs, even with bounded degree,\vadjust{\eject}
do not necessarily have the infinite collision
property.
This has been proved by Krishnapur and Peres \cite
{krishnapurperes2004}, looking at the comb lattice, that is
the subgraph of $\Z^2$ obtained by deleting all the horizontal edges
off the $x$-axis.

\subsection*{Infinite collision property of the percolation cluster}

We now prove that the infinite percolation cluster $C_{\infty} (\omega
)$ has the infinite collision property in $d = 2$,
which, by Lemma \ref{duality}, is equivalent to clustering of the
neutral competition model in two dimensions.
We use the same notation as before and let $\{X_t \}_t$ and $\{Y_t \}
_t$ be two independent continuous-time random walks run
at rate 1 on the infinite percolation cluster.
Let $W_t = (X_t, Y_t)$ and $\W_n = (\X_n, \Y_n)$ denote the
discrete-time Markov chain on
$C_{\infty} (\omega) \times C_{\infty} (\omega)$ with transition
probabilities
\begin{eqnarray*}
&&P \bigl(\W_{n + 1} = (A', B') \mid \W_n = (A, B)\bigr) \\
&&\qquad= \tfrac{1}{2} \bigl(q_1 (A, A')
\ind\{B = B' \} + q_1 (B, B')
\ind\{A = A' \} \bigr),
\end{eqnarray*}
where $q_n (A, B)$ denotes the $n$-step transition probability of the
lazy symmetric random walk on the infinite percolation cluster.
That is, at each time step, one of the two coordinates of $\W_n$ is
chosen at random with probability $1/2$.
This coordinate then moves according to the uniform distribution on the
neighbors or stands still, both with probability $1/2$, while
the other coordinate does not change.
Note that, at each step, with probability $1/2$, the process $\W_n$
does not move at all.
Note also that the processes $\{W_t \}_t$ and $\{\W_n \}_n$ can be
coupled in such a way that the sequences of states visited
by both processes are equal.
In particular, invoking in addition the Markov property and the
Borel--Cantelli lemma, to prove the infinite collision
property, it suffices to prove that
\[
P \bigl(\X_n = \Y_n \mbox{ for some } n \geq1 \mid \W_0 = (A, B)\bigr)
= 1.
\]
The first key to proving the infinite collision property of the cluster
is the following theorem, which is the analog of Theorem 1
in \cite{barlow2004}.
We state the result in the general $d$-dimensional case, though we only
deal with the two-dimensional case in the rest of this section.
\begin{theor}
\label{th:barlow}
Let $p > p_c$.
Then, there exist a subset $\Omega$ of the set of the realizations
with probability one and a collection of random variables
$\{S_A \}_{A \in\Z^d}$ such that the following holds:
\begin{enumerate}
\item We have $S_A (\omega) < \infty$ for each $\omega\in\Omega$
and $A \in C_{\infty} (\omega)$.
\item There are constants $c_1, c_2, c_3, c_4 > 0$ such that, for all
$A, B \in C_{\infty} (\omega)$,
%
%
\begin{eqnarray}
\label{eq:clust1}
q_n (A, B) & \geq& c_1 n^{- d/2} \exp(- c_2 |A - B|^2 /
n)\nonumber\\
&&\eqntext{\mbox{whenever } |A - B| \vee S_A (\omega) \leq n,} \\[-10pt]
\\[-10pt]
q_n (A, B) & \leq& c_3 n^{- d/2} \exp(- c_4 |A - B|^2 /
n)\nonumber\\
&&\eqntext{\mbox{whenever } S_A (\omega) \leq n.}
\end{eqnarray}
\end{enumerate}
\end{theor}

The proof of Theorem \ref{th:barlow} follows the lines of the proof of
its analog in \cite{barlow2004} and only differs
in two points: first, we consider a discrete-time lazy random walk
instead of a continuous-time random walk, and second,
processes under consideration evolve on the infinite percolation
cluster of site percolation instead of bond percolation.
To prove the sub-Gaussian upper estimate, the idea is to use a
discrete-time version of \cite{mathieuremy2004}, Theorem
1.1,
and the results of~\cite{coulhongrigoryanzucca2005}, Sections 5, 6 and
8, while the proof of the sub-Gaussian lower estimate
follows closely the strategy of \cite{barlow2004}.
Note that the choice of a lazy random walk is motivated by the fact
that one cannot expect the lower bound to hold for
any time $n$ for a standard simple random walk.
This is due to the fact that it has period 2.
In order to avoid unnecessary complications, we prefer to deal with an
aperiodic random walk.

In the sequel, to simplify notation, we write sums starting from
(or ending at) possibly noninteger real numbers,
but it is tacitly understood that one must consider their integer part.
To prove the infinite collision property, we define
\begin{eqnarray*}
F (n) & = & \sum_{j = 0}^n 2^{-n} \pmatrix{n \cr j}
q_j (A, X) q_{n - j} (B, X), \\
F_{\rho} (n) & = & \sum_{j = \rho n}^{(1 - \rho) n} 2^{-n} \pmatrix{n
\cr j} q_j (A, X) q_{n - j} (B, X),
\end{eqnarray*}
where $\rho\in(0, 1/2)$ and $A, B, X \in C_{\infty} (\omega)$.
\begin{lemma}
\label{lem:collisions}
Fix $\rho\in(0, 1/2)$, $A, B, X \in C_{\infty} (\omega)$ and $\ep> 0$.
Then $F (n) \leq(1 + \ep) F_{\rho} (n)$ when $n$ is sufficiently
large depending on $\rho$, $A$, $B$, $X$ and $\ep$.
\end{lemma}
\begin{pf}
By the Hoeffding inequality (\cite{hoeffding1963}, Theorem 1), we have
\begin{eqnarray*}
F (n) - F_{\rho} (n) &\leq&\sum_{j = 0}^{\rho n} 2^{-n} \pmatrix{n \cr
j} + \sum_{j = (1 - \rho) n}^n 2^{-n} \pmatrix{n \cr j} 
\\
&\leq&2 \exp\bigl(- 2 n (1/2 - \rho)^2\bigr).
\end{eqnarray*}
%
Taking $n$ such that
\begin{eqnarray*}
\rho n &\geq&\sqrt{n} \vee S_A (\omega) \vee S_B (\omega)\\[2pt]
&\geq&
|X - A| \vee|X - B| \vee S_A (\omega) \vee S_B (\omega),
\end{eqnarray*}
we may use the first inequality in (\ref{eq:clust1}).
Letting $\Phi$ denote the cumulative distribution function of the
standard normal, and also applying the central\vadjust{\eject} limit theorem,\vadjust{\goodbreak}
we obtain
%
%
\begin{eqnarray}
\label{eq:clust2}
F_{\rho} (n) & \geq&
\sum_{j = \rho n}^{(1 - \rho) n}\! 2^{-n} \pmatrix{n
\cr j} \frac{c_1}{j} \exp\biggl(\!- \frac{c_2 |X - A|^2}{j}
\!\biggr) \frac{c_1}{n - j} \exp\biggl(\!- \frac{c_2 |X - B|^2}{n - j}
\!\biggr)\hspace*{-25pt}
\nonumber\\ & \geq&
\sum_{k = \rho n}^{(1 - \rho) n} 2^{-n} \pmatrix{n
\cr j} \frac{c_1^2}{j (n - j)} \exp\biggl(- \frac{c_2 n^2}{j (n - j)} \biggr) \nonumber\\[-8pt]\\[-8pt]
&\geq&
(2 c_1 / n)^2 \exp\bigl(- c_2 / \bigl((1 - \rho) \rho\bigr)\bigr) \sum_{j = \rho n}^{(1 -
\rho) n} 2^{-n} \pmatrix{n \cr j}
\nonumber\\
& \geq&
C_5 n^{-2} \bigl(2 \Phi\bigl((1 - 2 \rho) \sqrt n\bigr) - 1\bigr)
\nonumber
\end{eqnarray}
for some $C_5 < \infty$.
To conclude, observe that
\[
F (n) = \biggl(1 + \frac{F (n) - F_{\rho} (n)}{F_{\rho} (n)}
\biggr) F_{\rho} (n),
\]
while the previous estimates (\ref{eq:clust2}) imply
\[
\lim_{n \to\infty} \frac{F (n) - F_{\rho} (n)}{F_{\rho} (n)} \leq\lim
_{n \to\infty} \frac{n^2}{C_5} \frac{\exp(- 2 n (1/2 -
\rho)^2)}{\Phi((1 - 2 \rho) \sqrt n) - 1/2} = 0.
\]
This completes the proof.
\end{pf}

With Theorem \ref{th:barlow} and Lemma \ref{lem:collisions}
in hand, we are now ready to prove that the
infinite percolation cluster has the infinite collision property in the
sense described above, that is, considering
continuous-time random walks run at a~constant rate, say 1.
Our proof relies in addition on an argument of Barlow, Peres and Sousi
\cite{barlowperessousi2009} who studied the number of
collisions of discrete-time random walks moving simultaneously at each
time step.
In order to understand the duality properties of the competition model,
we need, in contrast, to consider a pair
of random walks in which only one walk chosen uniformly at random can
move while the other walk stands still, thus
mimicking the evolution of a pair of independent continuous-time random walks.
\begin{theor}
\label{th:collisions}
Fix a realization $\omega$.
Then, for all $A, B \in C_{\infty} (\omega)$,
\[
P \bigl(\card\{n \dvtx\X_n = \Y_n \} = \infty\mid \W_0 = (A, B)\bigr) = 1.
\]
\end{theor}
\begin{pf}
Let $\gamma> 0$ to be chosen later, and define
\[
I_k = \sum_{n = k}^{k^2} \mathop{\sum_{|X - A| \vee|X - B| <
\sqrt n}}_
{S_X (\omega) \leq\gamma} I (X, n),
\]
where $I (X, n) = 1$ if there is a collision at time $n$ at site $X$,
and $= 0$ otherwise.
The first step is to find bounds for the first and second moments of
$I_k$ when~$k$ is large.

\textit{Lower bound}: $\E(I_k) \geq C \log k$ for some constant $C
> 0$ which does not depend on $A$, $B$ and
for all $k \geq k_1 (A, B)$.
First, we fix $\rho\in(0, 1/2)$ and observe that
\begin{eqnarray*}
\E^{A, B} I (X, n) &=& P^{A, B} (\X_n = \Y_n = X) = \sum_{j = 0}^{n}
2^{-n} \pmatrix{n \cr j} q_j (A,
X) q_{n - j} (B, X)
\\
&\geq& \sum_{j = \rho n}^{(1 - \rho) n} 2^{-n} \pmatrix{n
\cr j} q_j (A, X) q_{n - j} (B, X).
\end{eqnarray*}
In the previous sum, $j$ and $n - j$ are larger than $\rho n$.
Hence, for $n \geq k$, if
\[
|X - A| \vee|X - B| < \sqrt{\rho n} \quad\mbox{and}\quad k \geq
\rho^{-1} \bigl(S_A (\omega) \vee S_B (\omega)\bigr),
\]
then $j \wedge(n - j) \geq|X - A| \vee|X - B| \vee S_A (\omega)
\vee S_B (\omega)$ so Theorem \ref{th:barlow} implies
\begin{eqnarray*}
\E^{A, B} (I_k) & \geq& \sum_{n = k}^{k^2} \mathop
{\sum_{
|X - A| \vee|X - B| < \sqrt{\rho n}}}_{S_X (\omega) \leq\gamma}
\sum_{j = \rho n}^{(1 - \rho) n} 2^{-n} \pmatrix{n \cr j} \frac
{c_1^2}{j (n - j)} \\
&&\hspace*{120pt}{}\times\exp\biggl[- c_2 \biggl(\frac{\rho n}{j} + \frac{\rho n}{n - j}
\biggr) \biggr] \\
& \geq&
\sum_{n = k}^{k^2} \mathop{\sum_{|X - A| \vee|X -
B| < \sqrt{\rho n}}}_{S_X (\omega) \leq\gamma}
\sum_{j = \rho n}^{(1 - \rho) n} 2^{-n} \pmatrix{n \cr j} \biggl(\frac{2
c_1}{n} \biggr)^2 \exp\bigl(- c_2 / (1 - \rho)\bigr).
\end{eqnarray*}
This and the central limit theorem imply that, for $k$ large depending
on~$A, B$,
%
%
\begin{eqnarray}
\label{eq:clust3}
\E^{A, B} (I_k) &\geq&\exp\bigl(- c_2 / (1 - \rho)\bigr)
\nonumber\\
&&{}\times \sum_{n = k}^{k^2} (2 c_1 / n)^2 \card\bigl\{X \in
C_{\infty} (\omega) \dvtx\nonumber\\[-8pt]\\[-8pt]
&&\hspace*{93.1pt}|X - A| \vee|X - B| < \sqrt{\rho n},\nonumber\\
&&\hspace*{166.1pt} S_X
(\omega) \leq\gamma\bigr\}.\nonumber
\end{eqnarray}
Now, by the ergodic theorem,
\begin{eqnarray*}
&&\lim_{n \to\infty} \frac{\card\{X \in C_{\infty} (\omega) \dvtx|X -
A| \vee|X - B| <
\sqrt{\rho n}, S_X (\omega) \leq\gamma\}}
{\card\{X \in\Z^2 \dvtx|X - A| \vee|X - B| < \sqrt{\rho n} \}}
\\
&&\qquad= P \bigl(X \in C_{\infty} (\omega), S_X (\omega) \leq\gamma\bigr).
\end{eqnarray*}
In particular, there exists a constant $\delta> 0$ that only depends
on the percolation parameter $p$ such that for all
$\gamma$ and $k$ sufficiently large, we have
%
%
\begin{equation}
\label{eq:clust4}
\card\bigl\{X \in C_{\infty} (\omega) \dvtx|X - A| \vee|X - B| < \sqrt
{\rho n}, S_X (\omega) \leq\gamma\bigr\} \geq\delta\rho n
\end{equation}
for $n \geq k$.
By (\ref{eq:clust3}) and (\ref{eq:clust4}), there exists $k_1 (A, B)$
large such that
%
%
\begin{equation}
\label{eq:clust5}
\E^{A, B} (I_k) \geq c_1^2 \exp\bigl(- c_2 / (1 - \rho)\bigr) \sum_{n
= k}^{k^2} \frac{\delta\rho}{n} \geq C_6 (\log k^2 - \log k) = C_6 \log
k\hspace*{-30pt}
\end{equation}
for a suitable $C_6 > 0$ not depending on $A$, $B$, and all $k \geq k_1
(A, B)$.

\textit{Upper bound}: $\E(I_k^2) \leq C (\log k)^2$ for some
constant $C < \infty$ which does not depend on $A$, $B$
and for all $k \geq k_2 (A, B)$.
First, we observe that, for $l \geq n$,
\begin{eqnarray*}
\E^{A, B} (I (X, n) I (Y, l)) & = & P^{A, B} (\X_n = \Y_n = X,
\X_l = \Y_l = Y)
\\ & = &
\E^{A, B} I (X, n) \E^{X, X} (Y, l - n)
\end{eqnarray*}
from which it follows that
\begin{eqnarray*}
\E^{A, B} (I_k^2) &\leq&2 \sum_{n = k}^{k^2} \sum_{l =
n}^{k^2} \mathop{\sum_{|X - A| \vee|X - B| < \sqrt n}}_{S_X (\omega)
\leq
\gamma}
\mathop{\sum_{Y \dvtx|X - Y| < \sqrt{l - n}}}_{S_Y (\omega) \leq
\gamma}
\E^{A, B} I (X, n) \\
&&\hspace*{174pt}{}\times\E^{X, X} I (Y, l - n).
\end{eqnarray*}
Since $I (Y_1, l - n) I (Y_2, l - n) = 0$ whenever $Y_1 \neq Y_2$, we
also have
\[
\sum_{Y \dvtx|X - Y| < \sqrt{l - n}} \E^{X, X} I (Y, l - n) \ind\{
S_Y (\omega) \leq\gamma\} \leq1.
\]
Therefore, by applying Lemma \ref{lem:collisions} twice with $\ep=
1$, we deduce that there exists $\gamma$ large such that
for all $k$ sufficiently large
{\fontsize{10.5pt}{11pt}\selectfont{\begin{eqnarray*}
\hspace*{-4pt}&&\E^{A, B} (I_k^2) \\
\hspace*{-4pt}&&\qquad\leq 2 \sum_{n = k}^{k^2} \mathop{\sum_{|X - A| \vee
|X - B| < \sqrt n}}_{S_X (\omega) \leq
\gamma} 2 \Biggl( \sum_{j = \rho n}^{(1 - \rho) n} 2^{-n} \pmatrix{n \cr j}
q_j (A, X) q_{n - j} (B, X) \Biggr)\\
\hspace*{-4pt}&&\qquad\quad\hspace*{88.8pt}{}\times \Biggl(\frac{\gamma}{\rho} + 2 \sum_{l
= n + \gamma/ \rho}^{k^2} \mathop{\sum_{Y \dvtx|X - Y| < \sqrt{l -
n}}}_{S_Y (\omega) \leq
\gamma}
\sum_{i = \rho(l - n)}^{(1 - \rho) (l - n)} 2^{- (l - n)}\\
\hspace*{-4pt}&&\qquad\quad\hspace*{82pt}\hspace*{189pt}{}\times \pmatrix{l -
n \cr i} \\
\hspace*{-4pt}&&\qquad\quad\hspace*{205pt}{}\times q_i (X, Y) q_{l - n - i} (X, Y)
\Biggr).
\end{eqnarray*}}}
Observing that in the sums over $j$ and $i$ above, we have
\begin{eqnarray*}
j \wedge(n - j) & \geq& \rho n \geq\rho k \geq S_A (\omega)
\vee S_B (\omega), \\
i \wedge(l - n - i) & \geq& \rho(l - n) \geq\rho\gamma/ \rho
= \gamma\geq S_X (\omega)
\end{eqnarray*}
for all $k$ large depending on $A$, $B$, $X$, Theorem \ref{th:barlow}
implies that
\begin{eqnarray*}
\E^{A, B} (I_k^2) &\leq& 8 \sum_{n = k}^{k^2} 5 n
\Biggl( \sum_{j = \rho n}^{(1 - \rho) n} 2^{-n} \pmatrix{n \cr j} \frac
{c_3^2}{j (n - j)} \exp\bigl(- c_4 n \bigl(j^{-1} + (n - j)^{-1}\bigr) \bigr)\Biggr)
\\[-2pt]
&&\hspace*{21pt}{}\times \Biggl(\frac{\gamma}{\rho} + \sum_{l = n
+ \gamma/ \rho}^{k^2} \sum_{m = 0}^{\infty} \card\bigl\{Y \dvtx m \sqrt{l -
n} \leq|X - Y| \\[-2pt]
&&\hspace*{166.1pt}< (m + 1) \sqrt{l - n} \bigr\}
\\[-2pt]
&&\hspace*{21pt}{}\times \sum_{i = \rho(l - n)}^{(1 - \rho) (l - n)} 2^{- (l
- n)} \pmatrix{l - n \cr i} \frac{c_3^2}{i (l - n - i)}\\[-2pt]
&&\hspace*{86.3pt}{}\times \exp\bigl(- c_4 m^2
\bigl(i^{-1} + (l - n - i)^{-1}\bigr) \bigr)\Biggr).
\end{eqnarray*}
In particular, there exists $k_2 (A, B)$ large such that
%
%
\begin{eqnarray}
\label{eq:clust6}
\E^{A, B} (I_k^2) &\leq& 8 \sum_{n = k}^{k^2} \frac{5 c_3^2}{\rho(1 -
\rho) n} \nonumber\\
&&\hspace*{21.6pt}{}\times \Biggl(\frac{\gamma}{\rho} + \sum_{l = n + \gamma
/ \rho}^{k^2} \sum_{m = 0}^{\infty} \frac{5 (m + 2)^2 c_3^2}{\rho(1 -
\rho) (l - n)} \nonumber\\[-8pt]\\[-8pt]
&&\hspace*{116.8pt}{}\times \exp\bigl(- 4
c_4 (m + 1)^2\bigr) \Biggr)
\nonumber\\
&\leq& C_7 \sum_{n = k}^{k^2} \frac{1}{n} \Biggl(\frac{\gamma}{\rho} + C_8
\sum_{l = n + \gamma/ \rho}^{k^2} \frac{1}{l - n} \Biggr) \leq C_9 (\log
k)^2\nonumber
\end{eqnarray}
for suitable constants $C_7, C_8, C_9 < \infty$ not depending on $A$,
$B$ and all $k \geq k_2 (A, B)$.

Let $k (A, B) = k_1 (A, B) \vee k_2 (A, B)$.
By (\ref{eq:clust5}) and (\ref{eq:clust6}) and the Paley--Zygmund inequality,
\begin{eqnarray*}
&&P \bigl(I_k > (C_6 / 2) \log k \mid \W_0 = (A, B)\bigr)\\
&&\qquad \geq P \bigl(I_k > \E(I_k) / 2 \mid
\W_0 = (A, B)\bigr) \\
&&\qquad\geq(\E^{A, B} (I_k))^2 / 4 \E^{A, B} (I_k^2) \geq C_6^2 / 4 C_9
= c > 0
\end{eqnarray*}
for $k = k (A, B)$ and where, as $C_6$ and $C_9$, the constant $c > 0$
does not depend on the starting
points of the random walks.
Then, we define a sequence of stopping times and sites as follows:
we start at $n_0 = 0$ and $(A_0, B_0) = (A, B)$, and for all $j \geq1$
we define inductively
\[
n_j = n_{j - 1} + k (A_{j - 1}, B_{j - 1}) \quad\mbox{and}\quad
(A_j, B_j) = (\X_{n_j}, \Y_{n_j}).
\]
We say that there is a success at round $j \geq1$ when
\[
\card\{n \in[n_{j - 1}, n_j) \dvtx\X_n = \Y_n \} > (C_6 / 2)
\log k (A_{j - 1}, B_{j - 1})
\]
and observe that, at each round, the success probability is larger than
$c > 0$.
In particular, the probability mass function of the number of successes
up to round $j \geq1$ is stochastically larger than
a Binomial random variable with parameters $j$ and $c > 0$, from which
it follows that the ultimate number of successes, thus
the ultimate number of collisions, is almost surely infinite.
\end{pf}

As previously explained, clustering of the neutral competition model in
two dimensions follows from the combination
of Lemma \ref{duality} and Theorem \ref{th:collisions}.


\section{\texorpdfstring{Proof of Theorem \protect\ref{selection}}{Proof of Theorem 3}}
\label{sec:selection}

This section is devoted to the proof of Theorem~\ref{selection} and
the analysis of the competition model in the
presence of selection.
In the asymmetric case, the main difficulty arises from the fact that
the transition rates are no longer linear with respect
to the local frequencies of each symbiont type, which leads a
priori to a failure of duality techniques.
In particular, the invadability of type 1 is proved in the next
subsection by comparing directly the forward evolution of the
competition model with a~gambler's ruin model.
In contrast, under the additional assumption $\alpha_2 = 0$,
extinction of the symbionts of type 2 is established by invoking
duality techniques which are available for what we shall call a
threshold version of the competition model.
With this duality relationship in hands, the result follows as in the
previous section from random walk estimates.


\subsection*{Invasion of type 1}

In order to prove the first part of Theorem \ref{selection}, we first
let $\bar\xi_t (X)$ be the number of type 1
symbionts in the host at $X \in\Z^d$ and set
\begin{eqnarray*}
p_t (X) & = & \biggl({2d \alpha_1 \bar\xi_t (X) +
\beta_1 \sum_{Y \sim X} \bar\xi_t (Y)}\biggr)\\
&&{}\times\biggl(2d \alpha_1 \bar\xi
_t (X) +
2d \alpha_2 \bigl(N - \bar\xi_t (X)\bigr)\\
&&\hspace*{17.8pt}{} + \beta_1 \sum_{Y \sim X} \bar\xi_t (Y)
+ \beta_2 \sum_{Y \sim X} \bigl(N - \bar\xi_t (Y)\bigr)\biggr)^{-1},
\\
q_t (X) & = & \frac{2d \alpha_1 \bar\xi_t (X) +
\beta_1 \sum_{Y \sim X} \bar\xi_t (Y)}{2d N (\alpha_1 +
\beta_1)}.
\end{eqnarray*}
Observe that $p_t (X) = q_t (X)$ in the neutral case $\alpha_1 =
\alpha_2$ and $\beta_1 = \beta_2$, and that
\[
p_t (X), q_t (X) = \cases{
0 \quad \mbox{if and only if}\quad \bar\xi_t(X) = \bar\xi_t(Y) = 0, &\quad for all $Y
\sim
X$,\cr
1 \quad \mbox{if and only if} \quad \bar\xi_t(X) = \bar\xi_t(Y) = N, &\quad for all $Y
\sim X$.}
\]
Note also that, since $q_t (X)$ can take at most $(N + 1) (2d N + 1)$
different values,
\[
q^- := \inf\{q_t (X) \dvtx q_t (X) \in(0, 1) \} > 0
\]
and
\[
q^+ := \sup\{q_t (X) \dvtx q_t (X) \in(0, 1) \} < 1.
\]
Denote by $N_t$ the number of type 1 symbionts present in the system at
time~$t$.
If the number of symbionts of type 1 at time 0 is finite, then
\[
N_t \rightarrow\cases{
N_t + 1, &\quad  at rate $\displaystyle\sum_X \bigl(N - \bar\xi_t (X)\bigr) p_t (X)$,\vspace*{2pt}\cr
N_t - 1, &\quad at rate  $\displaystyle\sum_X \bar\xi_t (X) \bigl(1 - p_t (X)\bigr)$,}
\]
where the sum is over all $X \in\Z^d$ such that $p_t (X) \in(0, 1)$.
Now, we observe that in the neutral case when $\alpha_1 = \alpha_2$
and $\beta_1 = \beta_2$, the embedded Markov chain associated
to $\{N_t \}_t$ is the simple symmetric random walk on $\Z_+$ absorbed
at 0 [note that on each edge $(x, y)$ the rates of invasion
between $x$ and $y$ are symmetric], therefore the two rates above are
equal and
\begin{eqnarray*}
\frac{1}{N} \sum_X \bar\xi_t (X) &=& \frac{1}{N} \sum_X \bar\xi_t (X) \bigl(1
- p_t (X)\bigr)
+ \frac{1}{N} \sum_X \bar\xi_t (X) p_t (X)
\\
&=& \frac{1}{N} \sum_X \bigl(N - \bar\xi_t (X)\bigr) p_t (X)
+ \frac{1}{N} \sum_X \bar\xi_t (X) p_t (X) \\
&=& \sum_X p_t (X).
\end{eqnarray*}
This implies that for all $\alpha_1$ and $\beta_1$ such that $\alpha
_1 + \beta_1 \neq0$, and for all configurations
%
%
\begin{equation}
\label{eq:select1}
\frac{1}{N} \sum_X \bar\xi_t (X) = \sum_X \frac{2d \alpha_1 \bar\xi_t
(X) + \beta_1 \sum_{Y
\sim X} \bar\xi_t (Y)}{2d N (\alpha_1 + \beta_1)} = \sum_X q_t (X).
\end{equation}
Note that the expression of $q_t (X)$ depends neither on $\alpha_2$
nor on $\beta_2$ therefore the equation above does not hold in
the neutral case only.
It is always true regardless of the choice of the reproduction and
transmission rates.
Now, assume that $\alpha_1 \geq\alpha_2$ and $\beta_1 > \beta_2$.
We want to show that in this case the process $\{N_t\}_t$ has a
positive drift.
We say that:\vadjust{\goodbreak}
\begin{enumerate}
\item site $X \in\Z^d$ is bad at time $t$ when $q_t (X) \in(0, 1)$
and $\bar\xi_t (Y) = N$ for all $Y \sim X$ and
\item site $X \in\Z^d$ is good at time $t$ when $q_t (X) \in(0, 1)$
and $\bar\xi_t (Y) \neq N$ for some $Y \sim X$.
\end{enumerate}
Note that if $N_0$ is finite then at any time $t$ the sets of good and
bad sites are both finite.
The first ingredient to proving the result is to observe that for any
site either good or bad
%
%
\begin{eqnarray}
\label{eq:select2}
p_t (X) & = &
\biggl(2d \alpha_1 \bar\xi_t (X) + \beta_1 \sum
_{Y \sim X} \bar\xi_t (Y)\biggr)\nonumber\\
&&\hspace*{0pt}{}\times
\biggl(2dN \alpha_1 + 2dN \beta_1 - (\alpha_1 - \alpha_2) \bigl(N - \bar
\xi_t (X)\bigr)\nonumber\\[-8pt]\\[-8pt]
&&\hspace*{77.3pt}{} - (\beta_1 - \beta_2) \sum_{Y \sim X} \bigl(N - \bar\xi_t
(Y)\bigr)\biggr)^{-1} \nonumber\\
& \geq&
\frac{2d \alpha_1 \bar\xi_t (X) + \beta_1 \sum
_{Y \sim X} \bar\xi_t (Y)}{2dN (\alpha_1 + \beta_1)} = q_t (X),
\nonumber
\end{eqnarray}
while if we assume in addition that $X$ is a good site, then
%
%
\begin{eqnarray}
\label{eq:select3}
p_t (X) & \geq& \frac{2d \alpha_1 \bar\xi_t (X)
+ \beta_1 \sum_{Y \sim X} \bar\xi_t (Y)}{2dN (\alpha_1 + \beta
_1)} \biggl[1 - \frac{\beta_1 - \beta_2}{2dN (\alpha_1 + \beta_1)}
\biggr]^{-1} \nonumber\\[-8pt]\\[-8pt]
& \geq&
q_t (X) (1 - c)^{-1} \qquad\mbox{where } c =
\frac{\beta_1 - \beta_2}{2dN (\alpha_1 +
\beta_1)}.\nonumber
\end{eqnarray}
The second ingredient is to observe that, by definition of the lower
bound~$q^-$ and upper bound $q^+$, we
have $q^- q_t (X_1) \leq q^+ q_t (X_2)$ for all $X_1, X_2$ so
%
%
\begin{eqnarray}
\label{eq:select4}
&&q^- (1 - c)^{-1} q_t (X_1) + q^+ q_t (X_2) \nonumber\\
&&\qquad \leq q^- q_t
(X_1) + q^+ (1 - c)^{-1} q_t (X_2) ,\nonumber\\[-8pt]\\[-8pt]
&&\bigl(q^- (1 - c)^{-1} + q^+\bigr) \bigl(q_t (X_1) + q_t (X_2)\bigr) \nonumber\\
&&\qquad \leq (q^- +
q^+) \bigl(q_t (X_1) + (1 - c)^{-1} q_t
(X_2)\bigr).\nonumber
\end{eqnarray}
The third ingredient is to observe that if $X$ is bad then all sites $Y
\sim X$ are good, so the number of bad sites is
at most equal to the number of good sites.
In particular, letting
\[
B = \{X \dvtx X \mbox{ is bad} \}
\]
and
\[
G = \{X
\dvtx X \mbox{ is good} \},
\]
there exists a subset $G^* \subset G$ with $\card(G^*) = \card(B)$.
Then, combining (\ref{eq:select1})--(\ref{eq:select4}) gives
\begin{eqnarray*}
\sum_X p_t (X) & = &
\sum_{X\ \mathrm{bad}} p_t (X) + \sum
_{X \ \mathrm{good}} p_t (X) \\ & \geq&
\sum_{X \in B} q_t (X) + (1 - c)^{-1} \sum_{X
\in G^*} q_t (X) + (1 - c)^{-1} \sum_{X \in G \setminus G^*} q_t (X)
\\ & \geq&
\frac{q^- (1 - c)^{-1} + q^+}{q^- + q^+} \sum_{X
\in B \cup G^*} q_t (X) + (1 - c)^{-1} \sum_{X \in G
\setminus G^*} q_t (X)
\\ & \geq&
\frac{q^- (1 - c)^{-1} + q^+}{q^- + q^+} \sum_X q_t (X) = \frac{q^- (1
- c)^{-1} + q^+}{q^- + q^+} \frac
{1}{N} \sum_X \bar\xi_t (X).
\end{eqnarray*}
In particular, we have
\begin{eqnarray*}
\sum_X \bigl(N - \bar\xi_t (X)\bigr) p_t (X) & \geq&
\frac{q^- (1 - c)^{-1} + q^+}{q^- + q^+} \times\sum_X \bar\xi_t (X) -
\sum_X \bar\xi_t (X) p_t (X)
\\ & \geq&
\frac{q^- (1 - c)^{-1} + q^+}{q^- + q^+} \times\sum_X \bar\xi_t (X)
\bigl(1 - p_t (X)\bigr).
\end{eqnarray*}
By comparing the previous inequality with the transition rates of $\{
N_t \}_t$ and applying the gambler's ruin formula, we
can conclude that, starting with $K$ symbionts of type 1, we have
\[
P \Bigl(\lim_{t \to\infty} N_t = \infty\Bigr) \geq1 -
\biggl(\frac{q^- + q^+}{q^- (1 - c)^{-1} + q^+} \biggr)^K > 0.
\]
This establishes the first part of Theorem \ref{selection}.


\subsection*{Extinction of type 2}

We now prove that, under the extra assumption $\alpha_2 = 0$, type 1
outcompetes type 2, which is the second part
of Theorem \ref{selection}.
Letting $C_N$ denote the Cartesian product of the regular lattice $\Z
^d$ and $\K_N$ the expression of the Markov generator reduces to
\begin{eqnarray*}
L_2 f (\xi) & = &
\sum_{x \in C_N} \frac{\alpha_1 f_1 (x) + \beta_1
g_1 (x)}{\alpha_1 f_1 (x) + \beta_1 g_1 (x) + \beta_2 g_2
(x)} [f (\xi_{x, 1}) - f (\xi)]
\\ &&{} + \sum_{x \in C_N} \frac{\beta_2 g_2 (x)}{\alpha_1
f_1 (x) + \beta_1 g_1 (x) + \beta_2 g_2 (x)} [f (\xi_{x, 2}) - f (\xi)].
\end{eqnarray*}
Introduce $\beta= \beta_2$ and $\kappa= \beta_1 - \beta_2 > 0$.
Let $\gamma= \kappa\cdot(2d N)^{-1}$ and observe that
\[
g_1 (x) \neq0 \quad\mbox{implies that}\quad \kappa g_1 (x) \geq
\kappa\cdot(N \deg\pi(x))^{-1} = \kappa\cdot(2d N)^{-1} = \gamma
\]
from which it follows that
\begin{eqnarray*}
\frac{\alpha_1 f_1 (x) + \beta_1 g_1 (x)}{\alpha
_1 f_1 (x) + \beta_1 g_1 (x) + \beta_2 g_2 (x)} & = &
\frac{\alpha_1 f_1 (x) + \beta g_1 (x) + \kappa
g_1 (x)}{\alpha_1 f_1 (x) + \beta+ \kappa g_1 (x)} \geq \frac{\beta
g_1 (x) + \gamma}{\beta+ \gamma},
\\
\frac{\beta_2 g_2 (x)}{\alpha_1 f_1 (x) + \beta_1
g_1 (x) + \beta_2 g_2 (x)} & = &
\frac{\beta g_2 (x)}{\alpha_1 f_1 (x) + \beta+
\kappa g_1 (x)} \leq \frac{\beta g_2 (x)}{\beta+ \gamma}.
\end{eqnarray*}
In particular, it suffices to prove that type 1 symbionts outcompete
type 2 symbionts when the dynamics are described by the new
process $\{\zeta_t \}_t$ with Markov generator
\begin{eqnarray*}
L_3 f (\zeta) & = &
\sum_{x \in C_N} \frac{\beta g_1 (x) + \gamma
}{\beta+ \gamma} \ind\{g_1 (x) \neq0 \} [f (\zeta_{x, 1}) - f (\zeta)]
\\ &&{} + \sum_{x \in C_N} \biggl(\frac{\beta g_2 (x)}{\beta
+ \gamma} \ind\{g_1 (x) \neq0 \} + \ind\{g_1 (x) = 0 \} \biggr)
[f (\zeta_{x, 2}) - f (\zeta)],
\end{eqnarray*}
where as previously $\zeta_{x, i}$ is obtained from $\zeta$ by
assigning the value $i$ to vertex~$x$ and leaving the state of all the other
vertices unchanged.
The process can be constructed graphically from the three collections
of independent processes introduced in Table \ref{tab:harris} in the
%
%
\begin{table}[b]
\caption{Graphical representation}
\label{tab:harris}
\tabcolsep=-1pt
\begin{tabular*}{1.01\tablewidth}{@{\extracolsep{\fill}}lcc@{}}
\hline
\textbf{Notation} & \textbf{Description} & \textbf{Interpretation} \\
\hline
$T_n (x)$ & Poisson process with parameter 1 & Times of an update at
vertex $x$ \\
$U_n (x)$ & Uniform random variable on $(0, 1)$ & Determining the new
symbiont type \\
$V_n (x)$ & Uniform random variable on $\{1, 2, \ldots, 2d N\}$ &
Determining the new symbiont type \\
\hline
\end{tabular*}
\end{table}
following manner.
First, we equip the set $C_N$ with a total order relation.
Then, at time $T_n (z)$, we have the following alternative:
\begin{enumerate}
\item In the case when $U_n (z) < p := \beta(\beta+ \gamma)^{-1}$,
we select a neighboring vertex uniformly at random, say $y$, using in an
obvious manner the uniform random variable $V_n (z)$ and the total
order on the set $C_N$, and then draw an arrow from vertex $y$ to
vertex $z$.
\item In the case when $U_n (z) > p := \beta(\beta+ \gamma)^{-1}$,
we draw a set of $2dN$ arrows starting from each of the vertices adjacent
to the host at site $\pi(z)$ and pointing at vertex $z$.
\end{enumerate}
We call the events in rules 1 and 2 above jumping event and branching
event, respectively.
The type of vertex $z$ is updated at times $T_n (z)$, $n \geq1$, with
$z$ becoming of type 1 if at least one of the arrows that point
at vertex~$z$ originates from a type~1, and of type 2 otherwise.
In particular, if all the neighbors are of type 2 then the new type is
2 while if at least one neighbor is of type 1 then the new type
is chosen uniformly at random from the neighbors with probability $p$
or type 1 with probability $1 - p$, which produces the suitable
transition rates.
We say that there exists a path from space--time point $(y, T - s)$ to
point $(x, T)$ if there are sequences of times and vertices
\[
s_0 = T - s < s_1 < \cdots< s_{n + 1} = T
\quad\mbox{and}\quad
x_0 = y,\qquad x_1, \ldots, x_n = x,
\]
such that the following two conditions hold:
\begin{enumerate}
\item for $i = 1, 2, \ldots, n$, there is an arrow from $x_{i - 1}$ to
$x_i$ at time $s_i$ and
\item for $i = 0, 1, \ldots, n$, there is no arrow that points at the
segments $\{x_i \} \times(s_i, s_{i + 1})$.
\end{enumerate}
We define a set-valued process $\{\hat\zeta_s (x, T) \}_s$ by setting
for all $0 \leq s \leq T$
%
%
\begin{equation}
\label{eq:pseudo}
\hat\zeta_s (x, T) = \{y \in C_N \dvtx\mbox{there is a path from
$(y, T - s)$ to $(x, T)$} \}.
\end{equation}
Note that the process $\{\hat\zeta_s (x, T) \}_s$ consists of a
system of branching coalescing random walks in which particles independently
jump at rate $p$ and branch at rate $1 - p$.
We refer to the right-hand side of Figure \ref{fig:voter} on page
\pageref{fig:voter} for a~picture.
The introduction of the process (\ref{eq:pseudo}) is motivated by the
following lemma, which is somewhat reminiscent of the duality relationship
between the biased voter model and branching coalescing random walks.
\begin{lemma}
\label{pseudo-dual}
Assume that $\zeta_0 (z) \!=\! 1$ for some $z \in\hat\zeta_T (x, T)$.
Then, \mbox{$\zeta_T (x) \!=\! 1$}.
\end{lemma}
\begin{pf}
Let $z \in\hat\zeta_T (x, T)$ with $\zeta_0 (z) = 1$.
Then there is a unique path from $(z, 0)$ to $(x, T)$.
Using the same notation as in the definition of a path, we introduce
the jump process
\[
X_t = x_i \qquad\mbox{for all } s_i \leq t < s_{i + 1}
\quad\mbox{and}\quad X_T = x.
\]
From the construction of the process $\{\zeta_t \}_t$, we have $\zeta
_t (X_t) \!=\! 1$ for all $0 \!\leq\! t \!\leq\! T$.
The lemma follows immediately by applying the equation at time $t = T$
since $X_T = x$.
\end{pf}

The rest of the proof relies on standard random walk estimates
supplemented with a rescaling argument similar to the one described
in Section \ref{sec:invasion}.
In short, introducing the spatial regions
\[
B (X, K) = \{z \in C_N \dvtx\pi(z) \in X + (- K, K)^d \} \qquad\mbox
{for } X \in\Z^d \mbox{ and } K \in\Z_+,
\]
the next objective is to show that, for all $\ep> 0$ and $T = K^2$,
there exists $K$ large such that
%
%
\begin{eqnarray}
\label{eq:growth}
&&P \bigl(\zeta_T (x) = 1 \mbox{ for all } x \in B (X, 3K) \mid \zeta
_0 (x) = 1 \mbox{ for all } x \in B (X, K)\bigr)\nonumber\\[-8pt]\\[-8pt]
&&\qquad \geq1 - \ep.\nonumber
\end{eqnarray}
In view of Lemma \ref{pseudo-dual}, inequality (\ref{eq:growth})
follows directly from the following result.
\begin{lemma}
\label{growth}
Let $T = K^2$. Then, there exists $C_{10} < \infty$ and $\gamma_{10}
> 0$ such that
\[
P \bigl(\hat\zeta_T (x, T) \cap B (X, K) = \varnothing\bigr) \leq C_{10}
\exp(- \gamma_{10} K) \qquad\mbox{for all } x \in B (X, 3K).
\]
\end{lemma}
\begin{pf}
Let $x \in B (X, 3K)$.
The idea is to define a random walk $\{W_s \}_s$ embedded in the system
of branching coalescing random walks and connecting, with probability
close to one, the space--time point $(x, T)$ to a vertex in the ball $B
(X, K)$ at time 0.
The process starts at $W_0 = x$.
To define the dynamics, we also introduce the projection $\mathcal W_s
= \pi(W_s)$.
Then, the random walk jumps at each time $s$ such that $T - s = T_n
(W_{s-})$ when we have the following alternative:
\begin{enumerate}
\item Jumping event: when $U_n (W_s) < p$, there is an arrow from a
vertex, say $z$, to $W_{s-}$.
Then, the random walk jumps to vertex $z$, that is, we set $W_s = z$.
\item Branching event: when $U_n (W_s) > p$, there are $2dN$ arrows
that point at~$W_{s-}$.
Then, the random walk jumps to one of the tails, chosen randomly and
uniformly, that make the random walk's projection $\mathcal W_s$
closer to the center $X$ of the ball.
\end{enumerate}
Note that $W_s \in\hat\zeta_s (x, T)$ for all $s \in(0, T)$, as
desired. Introduce
\[
Y_s^i = |\pi_i (\mathcal W_s) - \pi_i (X)| \qquad\mbox{for } i
= 1, 2, \ldots, d,
\]
where $\pi_i$ is the projection on the $i$th axis in $\Z^d$.
Since each vertex has at most $dN$ neighbors which are closer to the
center of the target ball for a~total of $2dN$ neighbors, we have
\begin{eqnarray*}
\lim_{ h \to0} h^{-1} \cdot P (Y_{s + h}^i = Y_s^i + 1 \mid Y_s^i > 0) &
\leq&
p / (2d) =: r, \\
\lim_{ h \to0} h^{-1} \cdot P (Y_{s + h}^i = Y_s^i - 1 \mid Y_s^i > 0) &
\geq&
p / (2d) + (1 - p) / d =: l.
\end{eqnarray*}
Therefore, $\{Y_s^i \}_s$ is stochastically smaller than the random
walk $\{Z_s \}_s$ with a reflecting boundary at zero and
that otherwise jumps to the right at rate~$r$ and to the left at rate $l$.
Let $\tau$ denote the first time the random walk $Z_s$ hits the
boundary 0.
Since $r < l$, standard large deviation estimates for the Poisson
distribution imply that
%
%
\begin{equation}
\label{eq:fast}
P (\tau> C_{11} K \mid Z_0 \leq3K) \leq C_{12} \exp(-
\gamma_{12} K)
\end{equation}
for suitable constants.
Finally, we introduce the reverse asymmetric random walk $\{\bar Z_s \}
_s$ with state space $\Z$ that jumps to the right at rate $l$
and to the left at rate $r$. Letting
\[
u_k = P (\bar Z_s = 0 \mbox{ for some } s > 0 \mid \bar Z_0 =
k)
\]
a first-step analysis gives $l (u_{k + 1} - u_k) = r (u_k - u_{k -
1})$ and then
\[
(1 - a) (u_k - u_0) = (1 - a) \sum_{j = 0}^{k - 1} a^j
(u_1 - u_0) = (1 - a^k) (u_1 - u_0),
\]
where $a = r \cdot l^{-1}$.
It is straightforward to deduce that $u_k = a^k$.
In particular,
%
%
\begin{eqnarray}
\label{eq:slow}
&&P (Z_T > K \mid \tau< T) \leq P (Z_T > K \mid Z_0 = 0)
\nonumber\\[-8pt]\\[-8pt]
&&\qquad\leq
P (\bar Z_s = 0 \mbox{ for some } s < T \mid \bar Z_0 > K) \leq
u_K = a^K.
\nonumber
\end{eqnarray}
In conclusion, recalling the definition of the processes $\{Y_s^i \}
_s$, using the stochastic domination mentioned above, and
applying (\ref{eq:fast}) and (\ref{eq:slow}), we obtain
\begin{eqnarray*}
&&P \bigl(\hat\zeta_T (x, T) \cap B (X, K) = \varnothing\bigr) \\
&&\qquad\leq P \bigl(W_T \notin B
(X, K)\bigr) \leq d \times P (Y_T^i > K)
\\
&&\qquad\leq d \times P (\tau> T \mid Z_0 \leq3 K) + d \times P (Z_T > K
\mid \tau< T)
\\
&&\qquad\leq d \times C_{12} \exp(- \gamma_{12} K) + d \times a^K
\end{eqnarray*}
for all $K$ large.
Since $a = r \cdot l^{-1} < 1$, the lemma follows.
\end{pf}

Lemmas \ref{pseudo-dual} and \ref{growth} imply that, when viewed
under suitable scales, the set of space--time boxes which are void of
type 2
dominates oriented site percolation with parameter $1 - \varepsilon$.
This almost produces the second part of Theorem \ref{selection}.
The last problem is that oriented site percolation has a positive
density of unoccupied sites.
To prove that there is an in-all-directions expanding region which is
indeed void of type 2 symbionts, we apply a result
from Durrett \cite{durrett1992} which shows that unoccupied sites do
not percolate when $\varepsilon$ is close enough to 0.
Since symbionts of either type cannot appear spontaneously, once a
region is void of one type, this type can only reappear in the region
through invasion from the outside.
This then implies that our process has the desired property and
completes the proof of the second part of Theorem \ref{selection}.


\section*{Acknowledgments}
The authors would like to thank Martin Barlow, Yuval Peres and Perla
Sousi for fruitful discussions about random walks on random graphs,
as well as two anonymous referees for comments that helped to improve
the article.

%

%
\printaddresses

\end{document}